\documentclass{amsart} 
\usepackage{amsthm, amsfonts, amssymb, amsmath,epsfig,stmaryrd} 
\newtheorem{theorem}{Theorem}[section]
\newtheorem{lemma}[theorem]{Lemma}
\newtheorem{corollary}[theorem]{Corollary}
\newtheorem{proposition}[theorem]{Proposition}
\newtheorem{remark}[theorem]{Remark}

\theoremstyle{definition}
\newtheorem{definition}[theorem]{Definition}
\newtheorem{example}[theorem]{Example}
\newtheorem{remark/example}[theorem]{Remark/Example}

\newtheorem{question}[theorem]{Question}
\newtheorem{claim}[theorem]{Claim}

\numberwithin{equation}{section}

\def\NN{ {\bf N} }
\def\ZZ{ {\bf Z} }
\def\RR{ {\bf R} }

\def\PP{{\bf P}}
\def\AA{{\mathcal A}}
 
\newcommand{\ini}{\operatorname{in}}

\newcommand{\GL}{\operatorname{GL}}

\newcommand{\depth}{\operatorname{depth}}

\newcommand{\Tor}{\operatorname{Tor}}

\newcommand{\rank}{\operatorname{rank}}

\newcommand{\projdim}{\operatorname{projdim}}
\newcommand{\reg}{\operatorname{reg}}

\newcommand{\conv}{\operatorname{conv}}

\newcommand{\Flag}{\operatorname{Flag}}

\newcommand{\length }{\operatorname{length}}
\newcommand{\height }{\operatorname{height}}
\newcommand{\core }{\operatorname{core}}
\newcommand{\Pfaff }{\operatorname{Pfaff}}
\newcommand{\Plu }{\operatorname{Plu}}
\newcommand{\Segre }{\operatorname{Segre}}
\newcommand{\Ver }{\operatorname{Ver}}
\newcommand{\DetSym }{\operatorname{DetSym}}
\newcommand{\DetGen }{\operatorname{DetGen}}
\newcommand{\Simplex }{\operatorname{Simplex}}
\newcommand{\Grass }{\operatorname{Grass}}

\newcommand{\lt} {\lrcorner \!\!\!\!\shortuparrow}
\newcommand{\rt} {\Rsh}
\newcommand{\nlt}{\curvearrowleft}
\newcommand{\nrt}{\hookleftarrow}

\def\operator#1#2{\def#1{\mathop{\kern0pt\fam0#2}\nolimits}}  
\operator\length{length}
\operator\depth{depth}
\operator\ini{in}
\operator\reg{reg}
\operator\sind{sat}
\operator\height{height}
\operator\Tor{Tor}
\operator\projdim{proj.dim}
\operator\GL{GL}
\operator\ini{in}
\operator\gin{gin}
\operator\rl{rl}

\title [Nice Initial Ideals]
{Nice Initial Complexes of some Classical Ideals}

\author{Aldo Conca}
\address{Dipartimento di Matematica, 
Universite degli Studi di Genova,
Via Dodecaneso, 35
16146 Genova, Italy}
\email{conca@dima.unige.it}
\author{Serkan Ho\c{s}ten}
\address{Mathematics Department, San Francisco State University,
1600 Holloway Avenue, 
San Francisco, CA 94132, USA}
\email{serkan@math.sfsu.edu}
\author{Rekha R. Thomas}
\address{Department of Mathematics, University of
  Washington, 
Seattle, WA 98195,USA}
\email{thomas@math.washington.edu}
\date{\today}

\begin{document}

\begin{abstract} 
  This is a survey article on Gorenstein initial complexes of
  extensively studied ideals in commutative algebra and algebraic
  geometry. These include defining ideals of Segre and Veronese
  varieties, toric deformations of flag varieties known as Hibi
  ideals, determinantal ideals of generic matrices of indeterminates,
  and ideals generated by Pfaffians of generic skew symmetric
  matrices.  We give a summary of recent work on the construction of
  squarefree Gorenstein initial ideals of these ideals when the ideals
  are themselves Gorenstein. We also present our own independent
  results for the Segre, Veronese, and some determinantal cases.
\end{abstract}
\maketitle

\section{Introduction} \label{intro}  
Let $I$ be a homogeneous ideal in a polynomial ring $R$ over an
infinite field $K$ and let $\beta_{ij}(I)$ be the $(i,j)$-th Betti
number of $I$. Since passing to an initial ideal is a flat deformation
\cite[Chapter 15]{E}, $\beta_{ij}(I)\leq \beta_{ij}(\ini(I))$ for all
$i,j$ and every initial ideal $\ini(I)$ of $I$.  There are many
classes of toric or determinantal ideals arising from classical
constructions which are known to be minimally generated by Gr\"obner
bases in their special coordinate systems for carefully chosen term
orders.  These include the ideals defining Segre products and Veronese
subrings of polynomial rings, the ideals of minors of generic or
generic symmetric matrices of indeterminates, the ideals of Pfaffians
of generic skew symmetric matrices, defining ideals of Grassmannians
given by Pl\"ucker relations, etc. For any such classical ideal $I$
there is an explicit initial ideal, $\ini_{\mathrm{cla}}(I)$ (called
the {\em classical initial ideal} of $I$), which is squarefree and
Cohen-Macaulay, and has as many minimal generators as $I$ in each
degree, that is, $\beta_{0j}(I)=\beta_{0j}(\ini_{\mathrm{cla}}(I))$
for all $j$ (see \cite{BC,HT, H, Stu,Stu1}).  However, in most cases
$\beta_{ij}(I)\neq \beta_{ij}(\ini_{\mathrm{cla}}(I))$ for some $i$
and $j$.  In fact, the breakdown usually happens already at the first
syzygies; see Example \ref{spec} below. Therefore we are led to ask
the following question.

\begin{question}\label{q1}  Given a classical ideal $I$,  does there
  exist an initial ideal $\ini(I)$ such that
  $\beta_{ij}(I)=\beta_{ij}(\ini(I))$ for all $i,j$?
\end{question}  
 
\begin{example}\label{maxmin}  Question \ref{q1} has a positive answer
  for some instances, such as when $I$ is either the ideal of
  $m$-minors of a generic $m\times n$ matrix or the ideal of
  $(n-1)$-minors of a generic symmetric $n\times n$ matrix. In these
  cases, the classical initial ideal satisfies the conditions of
  Question \ref{q1}. The reason follows from a general remark. Let $J$
  be a homogeneous ideal in a polynomial ring $R$ which is generated
  by polynomials of degree $d$ and higher and has codimension $h$.  We
  denote the degree or multiplicity of $R/J$ by $\deg R/J$. If $R/J$
  is Cohen-Macaulay then it is easy to see that $\deg R/J \geq {h+d-1
  \choose d-1}$.  If $\deg R/J={h+d-1 \choose d-1}$ we say that $R/J$
  (or $J$) has minimal multiplicity (with respect to its initial
  degree).  Equivalently, a Cohen-Macaulay ring $R/J$ defined in
  degree $d$ and higher has minimal multiplicity if the quotient ring
  of $R/J$ by a regular sequence of $\dim R/J$ elements of degree one
  (its {\em Artinian reduction}) is isomorphic to
  $K[x_1,\dots,x_h]/\langle x_1,\dots,x_h \rangle^d$.  It follows that
  the Betti numbers of $J$ are equal to those of $\langle
  x_1,\dots,x_h \rangle^d$. In particular, if $J$ is Cohen-Macaulay of
  minimal multiplicity and $\ini(J)$ is a Cohen-Macaulay initial ideal
  then $\beta_{ij}(J)=\beta_{ij}(\ini(J))$ for all $i,j$. The ideal of
  $m$-minors of a $m\times n$ generic matrix and the ideal of
  $(n-1)$-minors of a generic symmetric $n\times n$ matrix are
  Cohen-Macaulay of minimal multiplicity and their classical initial
  ideals are Cohen-Macaulay. So the remark applies in these cases.
\end{example} 

The above examples do not represent the typical behavior.  Most
classical initial ideals do not have the correct Betti numbers. But
one can look for other non-classical initial ideals with the property
required in Question \ref{q1}.

\begin{example}\label{spec} Let $I$ be  the ideal of $2$-minors of the
  generic $3\times 3$ matrix $X=(x_{ij})$ and let 
$$\ini_{\mathrm{cla}}(I) = 
  \langle x_{11}x_{22},  x_{11}x_{23}, 
          x_{11}x_{32},  x_{11}x_{33},
          x_{12}x_{23},  x_{12}x_{33},
          x_{21}x_{32},  x_{21}x_{33},  
          x_{22}x_{33} \rangle
$$ 
be its classical initial ideal with respect to a diagonal term order. 
The Betti diagrams of $I$ and $\ini_{\mathrm{cla}}(I)$ are respectively: 
$$
\begin{array}{rrrrrrrrrr}  9& 16&  9& 0\\ 0&  0&  0& 1 \\ \\
\end{array} \quad \quad \mbox{and} \quad \quad 
\begin{array}{rrrrrrrrrr}  9& 16& 10& 2\\ 0&  1&  2& 1 \\ \\
\end{array}.
$$
Now we replace $\ini_{\mathrm{cla}}(I)$ with another initial ideal
$\ini_\succ(I)$ with respect to a reverse lexicographic term order
$\succ$ where the diagonal variables $x_{11}, x_{22}, x_{33}$ are 
smallest. The corresponding initial ideal is 
$$ 
\langle x_{23}x_{32}, x_{21}x_{32}, x_{13}x_{32}, x_{23}x_{31}, 
x_{13}x_{31}, x_{12}x_{31}, x_{12}x_{23}, x_{13}x_{21}, x_{12}x_{21} \rangle.
$$
One can check that the Betti diagram of $\ini_\succ(I)$ is identical
to that of $I$.

\end{example}  

In Section \ref{minors} we will generalize the phenomenon in Example
\ref{spec} to the ideal of $(n-1)$-minors in a generic $n\times n$
matrix. However, contrary to the above examples, in many classical
cases the answer to Question \ref{q1} is negative. In fact, the
property that is asked for may fail to hold for {\em all} initial
ideals (in the given coordinates).

\begin{example}\label{44cats}  Let $I$ be the ideal  of $2$-minors of
  a generic $4\times 4$ matrix. All initial ideals of $I$ are
  squarefree, Cohen-Macaulay, and generated in degree $\leq 4$, see
  \cite{Stu1}.  With the help of the the software package CaTS
  \cite{J} we computed all $4494288$ monomial initial ideals of $I$.
  They come in $4219$ distinct orbits modulo symmetries, and only
  $920$ orbits represent quadratically generated initial ideals.
  Computations in CoCoA \cite{Cocoa} reveal that the number of
  quadratic first syzygies of these initial ideals varies between $2$
  and $25$. Since $I$ has only linear first syzygies, that is
  $\beta_{1j}(I)=0$ for $j>2$, it follows that there is no initial
  ideal $\ini(I)$ such that $\beta_{ij}(\ini(I))=\beta_{ij}(I)$ for
  $i=0,1$ and all $j$.
\end{example}

The next best thing that one could ask for is an initial ideal which
has the correct number of generators and the correct Cohen-Macaulay
type. We also insist on asking for squarefree initial ideals so that
they can be represented by simplicial complexes.  Now we can state the
main question that this article addresses.

\begin{question}\label{fque}  Given a classical ideal $I$  which is
  Gorenstein does there exist a Gorenstein squarefree initial ideal
  of $I$ with the same number of generators as $I$?
\end{question}

Below is the list of classical ideals for which we study Question
\ref{fque}. The first three are examples of toric ideals which we
review now.  A polytope $P \subset \RR^d$ is called a {\em lattice
polytope} if its vertices lie in $\ZZ^{d}$.  Consider the embedding of
$P$ in $\RR^{d+1}$ given by $P \times \{1\} := \{ (p,1) \in \RR^{d+1}
\,:\,p \in P \}$ and let $C(P) \subset \RR^{d+1}$ be the cone over $P
\times \{1\}$. Then $M(P) := C(P) \cap \ZZ^{d+1}$ is a monoid whose
monoid algebra is $K[M(P)] := K[x^{m} \,:\, m \in M(P)]$, where $K$ is
an arbitrary field and $x=(x_1, \ldots, x_{d+1})$. The algebra
$K[M(P)]$ is graded by the exponent of $x_{d+1}$. Since $M(P)$ is
finitely generated as a monoid, $K[M(P)]$ is finitely generated as a
$K$-algebra. We say that $P$ is {\em normal} if $M(P)$ is generated by
the lattice points in $P \times \{1\}$ and hence $K[M(P)]$ is
generated by its monomials of degree one. A sufficient condition for
the normality of $P$ is the existence of a {\em unimodular
triangulation} of the lattice points in $P$.

Let $\mathcal P$ be the vector configuration consisting of the lattice
points in $P \times \{1\}$. Then the {\em toric ideal} of $\mathcal P$
is the homogeneous ideal $I_{\mathcal P} = \langle y^u - y^v \,:\,
\sum_{p_i \in \mathcal P} p_i u_i = \sum_{p_i \in \mathcal P} p_i v_i,
\, u_i, v_i \in \NN \rangle$ in the polynomial ring $K[y]$ where
$y = (y_1, \ldots, y_s)$ and $s = |\mathcal P|$. When $P$ is
normal, $I_{\mathcal P}$ is the presentation ideal of the algebra
$K[M(P)]$.  See \cite{Stu1} for details on toric ideals of vector
configurations.

If $K[M(P)]$ is Gorenstein, we say that $P$ is a {\em Gorenstein
  polytope}. Let $\textup{int}(M(P))$ denote the lattice points in the
interior of $C(P)$. It is well known that $K[M(P)]$ is Gorenstein if
and only of there exists $u \in \textup{int}(M(P))$ such that
$\textup{int}(M(P)) = u + M(P)$ \cite[Chapter 6]{BH}.

\begin{enumerate}
\item {\bf Segre($m,n$)}: Consider the Segre embedding of $\PP^{m-1}
  \times \PP^{n-1}$ in $\PP^{mn-1}$ parametrized by the monomial map
  $$K[x_{ij}] \rightarrow K[r_1, \ldots, r_m, s_1, \ldots, s_n], \,\,
  x_{ij} \mapsto r_is_j.$$ This is a toric variety with $P$ equal to
  the product of a standard $(m-1)$-dimensional simplex and and a
  standard $(n-1)$-dimensional simplex.  The corresponding vector
  configuration is $\mathcal P = \{e_i \oplus e'_j \,:\, 1 \leq i \leq
  m, \,\, 1 \leq j \leq n\}$ where $\{e_i\}$ and $\{e'_j\}$ are the
  standard unit vectors of $\RR^m$ and $\RR^n$ respectively.  Note
  that in this case, $P$ is a $(m+n-1)$-dimensional polytope that lies
  on the hyperplane $\sum r_i + \sum s_j = 2$ in $\RR^{m+n}$ and hence
  we can take $\mathcal P$ to be just the lattice points in $P$ as
  opposed to those in $P \times \{1\}$. The toric ideal $I_{\mathcal
  P}$ is generated by the $2$-minors of the $m \times n$ matrix
  $(x_{ij})$ of indeterminates. The polytope $P$ is Gorenstein with $u
  = (1,1,\ldots,1)$ if and only if $m=n$ \cite{BH,GW}. In this case, we
  will denote the defining ideal by $I(2,n)$.
  
\item {\bf Veronese($r,n$)}: Consider the $r$th Veronese embedding of
  $\PP^{n-1}$ in $\PP^N$ where $N = {r+n-1 \choose r-1}$ and $r \in
  \NN \backslash \{0,1\}$. This defines the toric ideal $I_{\mathcal
  P}$ where $P$ is the convex hull of all lattice points in $\NN^n$
  whose coordinates sum to $r$. The polytope $P$ is
  $(n-1)$-dimensional and lies on the hyperplane $\sum x_i = r$ in
  $\RR^n$.  The ideal $I_{\mathcal P}$ is Gorenstein if and only if
  $r$ divides $n$ \cite{GW}. When $r=2$, $I_{\mathcal P}$ is generated
  by the $2$-minors of a symmetric $n \times n$ matrix of
  indeterminates. We will denote this ideal by $J(2,n)$ throughout the
  article.

\item{\bf Hibi($m,n$)}: Let $e_{ij}$ be the unit vectors in $\NN^{m \times n}$,
 and let 
 $$\mathcal{P}_{m,n} = \{ e_{1a_1} + e_{2a_2} + \cdots + e_{ma_m} \, :
 \, 1 \leq a_1 < a_2 < \cdots < a_m \leq n \}.$$ The monoid algebra
 $K[M(P)]$ defined by the polytope $P$ that is the convex hull of the
 vectors in $\mathcal{P}_{m,n}$ is known as a Hibi ring.  These rings
 are obtained as certain toric deformations of the coordinate ring of
 $G(m,n)$, the Grassmannian of $m$-planes in $K^n$. The defining toric
 ideal $I_{m,n}$ is always Gorenstein \cite{BV}.  In Section
 \ref{flags} we will define general Hibi rings and discuss some
 results of Reiner and Welker \cite{RW}. Moreover we will describe in
 detail those Hibi rings obtained as Sagbi deformations of general
 flag varieties.
 
\item{\bf Plu($m,n$)}: Let $X = (x_{ij})$ be a $m \times n$ matrix of
  indeterminates. We denote by $[a_1, \ldots, a_m]$ the $m$-minor of
  $X$ with column indices $1 \leq a_1 < \ldots < a_m \leq n$. The
  algebra $K[[a_1, \ldots, a_m] \, : \, 1 \leq a_i \leq n] $ is the
  coordinate ring of $G(m,n)$. There are the well-known Pl\"ucker
  relations among these minors, see for instance \cite[Lemma
  7.2.3]{BH}.
  We let $K[x_\alpha \, : \, \alpha = (a_1, \ldots, a_m),\,\, 1 \leq
  a_1 < \cdots < a_m \leq n]$ be the polynomial ring in as many
  variables as the $m$-minors of $X$. We also define the $K$-algebra
  homomorphism $x_\alpha \mapsto [a_1, \ldots, a_m]$. The kernel
  $\Plu(m,n)$ of this map contains the quadratic polynomials that are
  preimages of the Pl\"ucker relations. Indeed, the preimages of the
  Pl\"ucker relations generate $\Plu(m,n)$.  This ideal is always
  Gorenstein \cite{Ful}.
  
\item{\bf DetGen($t,m,n$)}, {\bf DetSym($t,n$)}, and {\bf
    Pfaff($t,n$)}: Let $X$ be a $m \times n$ matrix of indeterminates.
    The ideal of all $t$-minors ($t>1$) of $X$ is called a
    determinantal ideal, and we will denote this ideal by
    $\DetGen(t,m,n)$. This ideal is Gorenstein if and only if $m= n$
    \cite{BV}.  Similarly, the ideal $\DetSym(t,n)$ will denote the
    ideal of $t$-minors of an $n\times n$ symmetric matrix $X$ of
    indeterminates.  This ideal is Gorenstein if and only if $n-t$ is
    even \cite{Goto}. Finally, for an even integer $t$ we let
    $\Pfaff(t,n)$ be the ideal of Pfaffians of order $t$ of a skew
    symmetric $n \times n$ matrix $X$.  The ideal of Pfaffians is
    always Gorenstein \cite{KL,Av}.
\end{enumerate}

This paper is organized as follows. In Section \ref{gen} we recall
general facts on Stanley-Reisner rings and Gorenstein simplicial
complexes.  Section \ref{gortoric} presents very recent results of
Bruns and R\"omer which imply that a Gorenstein toric ideal with a
squarefree initial ideal possesses a squarefree initial ideal that is
Gorenstein. All the toric ideals described above fall into this
category. However, the construction of Bruns and R\"omer does not
completely answer Question \ref{fque}. This is because one does not
know the degrees of the generators of the Gorenstein initial ideal
that exists via their result. We will treat the toric ideals $I(2,n)$
and $J(2,n)$ more extensively in Sections \ref{Segre} and
\ref{Veronese} and answer Question \ref{fque} positively in these
cases. In both cases we will show that the corresponding ideal has a
reverse lexicographic squarefree Gorenstein initial ideal where the
core of the associated simplicial complex is the boundary complex of a
simplicial polytope.  We will explicitly describe the facets and a
two-way shelling of these simplicial complexes. In Section \ref{flags}
we will examine Hibi rings more closely as the deformations of general
flag manifolds.  Section \ref{minors} will construct Gorenstein
initial ideals of $\DetGen(n-1,n)$, and Section \ref{jean-marie-pfaff}
will give similar constructions for $\Pfaff(t,n)$.

Before we go on, we would like to point out that except for the results in Section
\ref{Segre} and Section \ref{Veronese} which use shellings, a common theme to 
the results in this paper is the construction of a simplicial complex
$\Delta$ such that after the cone points of $\Delta$ are removed the remaining
complex is a simplical sphere. The first appearence of this kind of
result in our context is the {\em equatorial complex} construction 
of Reiner and Welker \cite{RW} (see Section \ref{jean-marie-pfaff}). 
Athanasiadis' result  \cite{A} that compressed polytopes (in particular,
the Birkhoff polytopes) with a {\em special simplex} 
have unimodular triangulations with an equatorial
complex (see Section \ref{gortoric}) was inspired by this result. 
Subsequently, Bruns and R\"omer \cite{BR} generalized Athanasiadis
result to all Gorenstein polytopes with a unimodular triangulation.

\section{Stanley-Reisner rings and Gorenstein complexes}\label{gen}

In this section we will recall briefly from \cite{BH, Sta} a few
important facts on Stanley-Reisner rings and Gorenstein simplicial
complexes.

Let $\Delta$ be a simplicial complex and let $K[\Delta]$ be its
Stanley-Reisner ring.  The dimension of a face $F \in \Delta$ is
$|F|-1$ and the dimension of $\Delta$ is the maximal dimension of its
facets.  We call $\Delta$ a pure complex if all its facets have the
same dimension.  We denote by ${\mathcal F}(\Delta)$ the set of facets
of $\Delta$. Every simplicial complex has an (essentially unique)
geometric realization.  A simplicial complex $\Delta$ of dimension
$d-1$ is said to be a simplicial sphere if its geometric realization
is homeomorphic to the sphere $S^{d-1}\subset \RR^{d}$.  The Hilbert
series of $K[\Delta]$ where $\Delta$ is $(d-1)$-dimensional has the
form
$$ \frac{h_0 + h_1t + \cdots + h_st^s}{(1-t)^d}$$ with $h_i\in \ZZ$,
$h_s\neq 0$ and $s\leq d$. The vector $h(\Delta) := (h_0, h_1, \ldots,
h_s)$ is called the $h$-{\em vector} of $\Delta$.  The $a$-invariant
$a(K[\Delta])$ of $K[\Delta]$ is $s-d$, the degree of the Hilbert
series as a rational function.

Given subsets $F_1,\dots,F_k$ of a given set $V$ we denote by $\langle
F_1,\dots,F_k\rangle$ the smallest simplicial complex containing the
$F_i$. Furthermore, if $\Delta_1$ and $\Delta_2$ are simplicial
complexes on disjoint sets of vertices $V_1$ and $V_2$, the join of
$\Delta_1$ and $\Delta_2$ is $\Delta_1*\Delta_2 = \{ A \cup B \, : \,
A \in \Delta_1, \, B \in \Delta_2\}$.  We let $CP(\Delta) = \{v \in V
\,: \, v \in F \,\, \forall F \in {\mathcal F}(\Delta)\}$ be the
cone-points of $\Delta$.  We also let $\core(\Delta)$ be the
restriction of $\Delta$ to the set of vertices not in $CP(\Delta)$
. This implies that $\Delta=\core(\Delta)*\Simplex(CP(\Delta))$. Note
that the elements in $CP(\Delta)$ correspond exactly to those
variables which do not appear in the generators of the Stanley-Reisner
ideal of $\Delta$.  So $K[\Delta]$ is just a polynomial extension of
$K[\core(\Delta)]$.

A simplicial complex $\Delta$ is said to be Cohen-Macaulay or
Gorenstein with respect to the field $K$ if the Stanley-Reisner ring
$K[\Delta]$ is Cohen-Macaulay or Gorenstein. The {\em link} of a face
$F \in \Delta$ is $\textup{lk}_\Delta(F) = \{ G \in \Delta \,:\, G
\cup F \in \Delta, \, G \cap F = \emptyset \}$. 

\begin{theorem} \cite[Corollary 4.2]{Sta}, \cite[Theorem 5.1]{Sta}
A simplicial complex $\Delta$ is 
\begin{itemize}
\item Cohen-Macaulay over $K$ if and only if for all $F \in \Delta$
and all $i < \textup{dim}(\textup{lk}_\Delta(F))$, we have
$\tilde{H_i}(\textup{lk}_\Delta(F); K) = 0$, and 
\item Gorenstein over $K$
if and only if for all $F \in \core(\Delta)$, \\
$\tilde{H_i}(\textup{lk}_{\core (\Delta)}(F);K) = \left \{
\begin{array}{ll} K \,\, \mbox{if} \,\, i =
\textup{dim}(\textup{lk}_{\core (\Delta)}(F))\\ 0 \,\, \mbox{if} \,\, i
< \textup{dim}(\textup{lk}_{\core (\Delta)}(F))
\end{array} \right.$
\end{itemize}
\end{theorem}

A simplicial complex $\Delta$ of dimension $d-1$ is said to be
shellable if it is pure and ${\mathcal F}(\Delta)$ can be totally
ordered so that for every non-minimal $F\in \mathcal F(\Delta)$ the
simplicial complex
$$\langle F\rangle \cap \langle G \in\mathcal F(\Delta) : G<F\rangle
\eqno{(1)}$$ is pure of dimension $d-2$. The total order $<$ is called
a shelling of $\Delta$.  Shellable simplicial complexes are
Cohen-Macaulay (over any field) and their Hilbert series can be
described in terms of the facets of the simplicial complex (1) as $F$
varies. Important features of Gorenstein simplicial complexes are
summarized below.

\begin{lemma}\label{lem1} Let $\Delta$ be a simplicial complex.
\begin{itemize}
\item[(a)] If $K[\Delta]$ is Gorenstein, then
  $a(K[\Delta])=-|CP(\Delta)|$. Equivalently, if $K[\Delta]$ is
  Gorenstein and $\core(\Delta)$ has dimension $d-1$ then
  $h_d(\Delta)=1$ and $h_i(\Delta)=0$ for $i>d$.
\item[(b)]  If $\Delta$  is a simplicial sphere  then $K[\Delta]$ is
Gorenstein.
\end{itemize}
Furthermore assume that $\core(\Delta)$ is shellable and every face of
codimension $1$ (i.e.  dimension $\dim \core(\Delta)-1$) is contained
in exactly two facets. Then
\begin{itemize}
\item[(c)]   $\core(\Delta)$  is a simplicial sphere.
\item[(d)] $K[\Delta]$ is Gorenstein.
\end{itemize}
\end{lemma}

\begin{proof} 
For (a) and (b) see \cite[Section 5.6]{BH}.  Statement (c) is proved 
in \cite[4.7.22]{BLSWZ} and (d) follows from (b) and (c).
\end{proof}

A shelling $<$ of $\Delta$ is said to be a two-way shelling if the
facets in the reversed order also give a shelling.  Line shelling of
simplicial polytopes are typical examples of two-way shellings.  The
shellings that we describe in this paper are shown to be two-way (but
we do not know whether they are line-shellings). 

\section{Gorenstein Toric Ideals} \label{gortoric}

In this section we survey recent results on Gorenstein toric ideals
that are relevant to this paper.  We use the notation introduced
earlier on polytopes, monoid rings, and toric ideals.  The following
theorem was proved by Bruns and R\"omer \cite{BR} and relates to
Question \ref{fque} addressed in this paper.

\begin{theorem} \cite[Corollary 7]{BR} \label{BRcor7}
  Let $P$ be a Gorenstein lattice polytope such that the set of
  lattice points $\mathcal P$ in $P$ admit a regular unimodular
  triangulation. Then the toric ideal $I_{\mathcal P}$ has a
  squarefree Gorenstein initial ideal.
\end{theorem}

We summarize the key ideas in the proof of this theorem from
\cite{BR}.  If $\ini_\succ(I)$ is a monomial initial ideal of the
ideal $I$, then the radical ideal $\textup{rad}(\ini_\succ(I))$ is
squarefree and is the Stanley-Reisner ideal of a simplicial complex
$\Delta(\ini_\succ(I))$. The simplicial complex
$\Delta(\ini_\succ(I))$ is called the {\em initial complex} of $I$
with respect to $\succ$.  Theorem 8.3 in \cite{Stu1} proves that if
$\ini_\succ(I_{\mathcal P})$ is a monomial initial ideal of the toric
ideal $I_{\mathcal P}$, then the initial complex
$\Delta(\ini_\succ(I_{\mathcal P}))$ is precisely the {\em regular
triangulation} $\Delta_\succ(\mathcal P)$ of $\mathcal P$ induced by
$\succ$. Further, Corollary~8.9 in \cite{Stu1} shows that such an
initial ideal is squarefree if and only if $\Delta_\succ(\mathcal P)$
is unimodular. Thus when $\ini_\succ(I_{\mathcal P})$ is squarefree,
the ring $K[y]/\ini_\succ(I_{\mathcal P})$ is the
Stanley-Reisner ring $K[\Delta_\succ(\mathcal P]$.

The main theorem (Theorem 3) in \cite{BR} states that whenever $M$ is
a normal affine monoid such that the monoid algebra $K[M]$ is
positively graded and Gorenstein, then there exists a simpler normal 
affine Gorenstein monoid algebra $K[N]$ where the monoid $N$ is
obtained as a projection of $M$. In the case where $M = M(P)$ and $P$
is a normal Gorenstein lattice polytope, the monoid algebra $K[N]$ is
generated in degree one in the grading inherited from $K[M]$ and hence
equals $K[M(Q)]$ where $Q$ is the polytope spanned by the exponents of
the monomials in $K[N]$ of degree one. Further, $Q$ is a Gorenstein
lattice polytope with a unique interior lattice point. If we let  
the Hilbert series of $K[M(P)]$ be
$$ \frac{h_0 + h_1t + \cdots + h_dt^d}{(1-t)^{\textup{dim}(P)+1}}$$
and the $h$-vector $h(P) := (h_0, h_1, \ldots, h_d)$, then they show
that $h(P) = h(Q) = h(\partial(Q))$ where $\partial(Q)$ is the
boundary of $Q$ (see \cite[Corollary 4]{BR}).

In Theorem~\ref{BRcor7}, we are given a Gorenstein lattice polytope
$P$ such that $\mathcal P$, the lattice points of $P$, admits a
regular unimodular triangulation. Let $\Delta(\mathcal P)$ be the
induced regular unimodular triangulation of $\mathcal P$ and $J$ the
squarefree monomial initial ideal of $I_{\mathcal P}$ whose initial
complex is $\Delta(\mathcal P)$. Since $\mathcal P$ has a unimodular
triangulation, $P$ is normal and the polytope $Q$ constructed above
exists. Since the triangulation $\Delta(\mathcal P)$ is a regular
unimodular triangulation of $M(P)$, equivalently of $C(P)$, with all
cones generated by elements of $\mathcal P$, $M(Q)$ and hence $Q$
inherits a regular unimodular triangulation $\Delta(Q)$.  Project the
vertices of $\Delta(Q)$ on a sphere around the unique lattice point in
$Q$ and let $P'$ be the simplicial polytope obtained as the convex
hull of these projected vertices. Since $M(P)$ is Gorenstein there
exists a unique $x \in \textup{int}(M(P))$ such that
$\textup{int}(M(P)) = x + M(P)$. Let $p_1, \ldots, p_m$ be a subset of
the minimal generating set (Hilbert basis) of $M(P)$ such that $x =
p_1 + \ldots + p_m$. From the construction of $Q$ and $P'$ it follows
that $\Delta(\mathcal P)$ is the join of the simplicial complex
$\Delta(P')$ corresponding to $\partial(P')$ and the simplex with
vertices $p_1, \ldots, p_m$. This implies that the variables $y_1,
\ldots, y_m$ in $K[y]$ corresponding to $p_1, \ldots, p_m$ form a
regular sequence modulo $J$ and hence $K[y]/J$ is Gorenstein since
$\Delta(P')$ is the boundary of a simplicial polytope. We note that
the main goal of \cite{BR} was to prove that if $P$ is an integer
Gorenstein polytope whose lattice points admit a unimodular
triangulation then the $h$-vector of $P$ is unimodal.

We now apply Theorem~\ref{BRcor7} to various Gorenstein lattice
polytopes and their toric ideals listed in the Introduction.

\begin{enumerate}
\item {\bf Segre($n,n$)}: All regular triangulations of $\mathcal P$
  are known to be unimodular and $P$ is Gorenstein. Hence by
  Theorem~\ref{BRcor7}, $I_{\mathcal P}$ has a squarefree Gorenstein
  initial ideal. In Section~\ref{Segre} we will construct an explicit
  term order $\succ$ such that the initial ideal
  $\ini_\succ(I_{\mathcal P})$ is {\em quadratic}, squarefree and
  Gorenstein.
  
\item {\bf Veronese($r,n$)}: The polytope $P$ defining $I_{\mathcal
   P}$ is a simplex that admits a unimodular triangulation consisting
   of empty simplices whose facets are parallel to the facets of
   $P$. When $P$ is Gorenstein ($r$ divides $n$), Theorem~\ref{BRcor7}
   applies.  In Section \ref{Veronese}, in the case of $r=2$ and
   $n=2m$ we will exhibit an explicit initial ideal
   $\ini_\succ(J(2,n))$ that is quadratic, squarefree and Gorenstein.

\item{\bf Hibi($m,n$)}: The vector configuration ${\mathcal P}$
   defining the Hibi ring is affinely isomorphic to the vertices of
   the order polytope of the lattice of order ideals of the product of
   the chains $[m] \times [n-m]$, see Section~\ref{flags} or
   \cite[Remark 11.11]{Stu1}. Order polytopes are known to have
   unimodular triangulations, and since $I_{\mathcal P}$ is a toric
   deformation of $\Plu(m,n)$ the polytope $P$ is also Gorenstein.
   Again, Theorem~\ref{BRcor7} applies.

\end{enumerate}

There is one more polytope we have not mentioned so far which played a
motivating role for both Theorem~\ref{BRcor7} and the earlier work of
Athanasiadis \cite{A}. \\ 

\noindent {\bf Birkhoff($n$)}: Recall that the $n$th Birkhoff polytope
in $\RR^{n \times n}$ is the convex hull of all the $n \times n$
permutation matrices. In this case, $\mathcal P$ equals the set of
$n!$ vertices of this polytope. Birkhoff polytopes are known to be
{\em compressed} which means that all their reverse lexicographic
triangulations are unimodular \cite{Stan80}.  Further, they are also
Gorenstein. Hence again by Theorem~\ref{BRcor7}, $I_{\mathcal P}$ has
a squarefree Gorenstein initial ideal. In the rest of this section we
briefly describe Athanasiadis' method.

A {\em special simplex} of a lattice polytope $P \subset \RR^d$, is a
collection of vertices $\Sigma = \{v_1, \ldots, v_q\}$ of $P$ with the
property that every facet of $P$ contains all but one vertex in
$\Sigma$. For the $n$th Birkhoff polytope, the collection of
permutation matrices corresponding to the cyclic subgroup of $S_n$
generated by the cycle $(1 \, 2 \, 3 \, \cdots \, n)$ forms a special
simplex. To see this note that the facets of the $n$th Birkhoff
polytope are cut out by the hyperplanes $x_{ij} = 0$ in $\RR^{n \times
n}$ and each facet misses exactly one permutation in the above cyclic
group. Note that special simplices are not contained in the boundary
of $P$. If $V$ is a linear subspace in $\RR^d$, let $P/V$ denote the
quotient polytope equal to the image of $P$ under the canonical
projection $\RR^d \rightarrow \RR^d/V$.

\begin{lemma} \cite[Proposition 2.3]{A} \label{prop2.3Athan}
Let $P$ be a $d$-dimensional polytope in $\RR^d$ with a special
simplex $\Sigma$ such that $\mathcal P$ has a triangulation isomorphic
to $\Sigma \ast \Delta$. Let $V$ be the linear subspace parallel to
the affine span of $\Sigma$. Then the boundary complex of the quotient
polytope $P/V$ inherits a triangulation abstractly isomorphic to
$\Delta$ and its faces are precisely the faces of $P$ that do not
intersect $\Sigma$.
\end{lemma}

\begin{lemma} \cite[Lemma 3.4]{A} \label{lemma3.4Athan}
  Suppose that $v_1 \prec \cdots \prec v_q \prec \cdots \prec v_{p-1}
  \prec v_p$ is an ordering of the vertices of a lattice polytope $P$
  such that $\Sigma = \{ v_1, \ldots, v_q\}$ is a special simplex of
  $P$.  Let $\Delta$ be the reverse lexicographic triangulation of $\{v_p,
  \cdots, v_{q+1} \}$ with respect to the order $\succ$. Then
\begin{enumerate} 
\item The reverse lexicographic triangulation $\Delta_\succ(P)$ is
  isomorphic to $\Sigma \ast \Delta$, and 
\item $\Delta$ is isomorphic to the boundary complex of $P/V$ which in
  turn is isomorphic to the boundary complex of a simplicial polytope
  of the same dimension as $P/V$.
\end{enumerate}
\end{lemma}

We state a modified version of the main theorem in \cite{A}.

\begin{theorem} \cite[Theorem 3.5]{A} \label{thm3.5Athan}
  Suppose $P$ is a lattice polytope and 
$v_1 \prec \cdots \prec v_q \prec \cdots \prec v_{p-1}
  \prec v_p$ is an ordering of its
  vertices such that (i) $P$ is compressed and (ii) $\Sigma = \{ v_1,
  \ldots, v_q\}$ is a special simplex of $P$. Then the $h$-vector
  $h(P)$ equals $h(\partial(Q))$ where $Q$ is a simplicial polytope
  whose boundary is isomorphic to the reverse lexicographic
  triangulation of $\{v_p, \cdots, v_{q+1} \}$ with respect to the
  order $\succ$.
\end{theorem}

\begin{proof}
First we invoke the fact that if $\Delta$ is any unimodular
triangulation of $\mathcal P$, then $h(P) = h(\Delta)$. In the
situation of the theorem, since $P$ is compressed, $\Delta_\succ(P)$
is unimodular and hence $h(P) = h(\Delta_\succ(P))$. By
Lemma~\ref{lemma3.4Athan} (i), $\Delta_\succ(P) = \Sigma \ast \Delta$
where $\Delta$ is the reverse lexicographic triangulation of
$\{v_p,\cdots, v_{q+1} \}$ with respect to the order $\succ$. Thus
$$h(P) = h(\Delta_\succ(P)) = h(\Sigma \ast \Delta) = h(\Delta)$$
where the third equality is a standard fact
about joins of simplicial complexes and the last equality follows
since the $h$-vector of a simplex is always $1$. By
Lemma~\ref{lemma3.4Athan} (ii), $\Delta$ is isomorphic to the boundary
complex of a simplicial polytope $Q$ whose boundary is isomorphic to
the reverse lexicographic triangulation of $\{v_p, \cdots, v_{q+1} \}$
with respect to the order $\succ$ which completes the proof. 
\end{proof}

The following is a modified version of Corollaries 4.1 and 4.2 in
\cite{A} adapted to this paper. 

\begin{corollary} \label{cors4.1and2Athan} 
Let $P$ be a compressed Gorenstein lattice polytope such that $M(P)$
is generated by the vertices of $P$. Then the toric ideal $I_{\mathcal
  P}$ has a squarefree Gorenstein initial ideal. In particular, the 
toric ideal of the $n$th Birkhoff polytope has a squarefree
Gorenstein initial ideal. 
\end{corollary} 

\begin{proof} 
  Since $P$ is Gorenstein, there exists unique $x \in
  \textup{int}(M(P))$ such that $\textup{int}(M(P))= x + M(P)$. Let
  $v_1, \ldots, v_q$ be vertices of $P$ such that $x = v_1 + \ldots +
  v_q$. Athanasiadis proves that $\Sigma = \{v_1, \ldots, v_q\}$ is a
  special simplex of $P$ \cite[Corollary 4.1]{A}. Now consider any
  reverse lexicographic ordering of the vertices of $P$ such that $v_q
  \succ \ldots \succ v_1$ comes last in the ordering. Then the
  conclusion of Theorem~\ref{thm3.5Athan} holds. Let $J$ be the
  initial ideal $\ini_\succ(I_{\mathcal P})$. Since $\Delta_\succ(P) =
  \Sigma \ast \Delta$ (from Theorem~\ref{thm3.5Athan}) is the initial
  complex of $J$, $J$ is squarefree. Further, since $\Delta$ is the
  boundary complex of a simplicial polytope and $\Sigma$ is a simplex,
  $K[y]/J$ is Gorenstein.
\end{proof} 

Note that Theorem~\ref{BRcor7} is a generalization of
Corollary~\ref{cors4.1and2Athan}. Further examples of Gorenstein
lattice polytopes that satisfy the conditions of
Corollary~\ref{cors4.1and2Athan} can be found in \cite{HO2005}. We
will see in Section \ref{Segre} that the polytope $P$ of Segre($n,n$)
also satisfies the conditions of Corollary~\ref{cors4.1and2Athan}
providing yet another proof that its toric ideal has a squarefree
Gorenstein initial ideal.

\section{Gorenstein Segre products}\label{Segre}

As we indicated already, $I(2,n)$ is generated by the $2$-minors of a
$n\times n$ matrix $X=(x_{ij})$ of indeterminates, and it is  
an ideal of the polynomial ring $K[x_{ij}]$.
The Hilbert series of $K[x_{ij}]/I(2,n)$ is given by
$$\sum_i { n-1 \choose i}^2 z^i/(1-z)^{2n-1}.$$ So the $a$-invariant
is $n-1-(2n-1)=-n$, and therefore any squarefree Gorenstein initial
complex of $I(2,n)$ must have exactly $n$ cone points. The classical
initial ideal of $I(2,n)$ is the one associated to a ``diagonal" term
order, namely a term order which selects main diagonals as initial
terms of minors and it is generated by the products $x_{ij}x_{hk}$
with $i<h$ and $j<k$. The facets of this initial complex are the paths
from $(n,1)$ to $(1,n)$ in an $n \times n$ grid. Table \ref{44path}
shows a typical facet of the classical initial complex of
$I(2,4)$. Since $(n,1)$ and $(1,n)$ (corresponding to the variables
$x_{n1}$ and $x_{1n}$) are the only points that belong to every facet,
this initial complex has only two cone points. So for $n>2$ it is not
Gorenstein.
\begin{scriptsize}
\begin{table}[ht] 
\caption{}
\renewcommand\arraystretch{1.5}
\noindent
\[ 
\begin{array}{|c|c|c|c|c|}
\hline    &    &     & *   \\
\hline    &  * &  *  & *    \\
\hline    &  * &    &       \\
\hline  * &  * &    &       \\
\hline
\end{array} 
\] 
\label{44path}
\end{table}
\end{scriptsize}
In order to construct a Gorenstein initial complex we consider a term
order where the initial term of a minor is its main diagonal unless
the main diagonal of the minor involves elements of the main diagonal
of the matrix.  Formally, for every $i<h$ and $j<k$ the initial term
of the minor $x_{ij}x_{hk}-x_{ik}x_{hj}$ is $x_{ij}x_{hk}$ unless
$i=j$ or $h=k$.  We can define such a term order by a reverse
lexicographic order $x_{11} \prec x_{22} \prec \cdots \prec x_{nn}
\prec \{ x_{ij} : i \neq j\}$ where the latter set of variables are
ordered so that $x_{ij}\succ x_{hk}$ if $|i-j|<|h-k|$.  For instance
for $n=4$, we could use: $ x_{12}\succ x_{21}\succ x_{23}\succ
x_{32}\succ x_{34}\succ x_{43}\succ x_{13}\succ x_{24}\succ x_{31}
\succ x_{42}\succ x_{14}\succ x_{41}\succ x_{44}\succ x_{33}\succ
x_{22} \succ x_{11}.  $

The initial terms of the $2$-minors are the monomials in the variables
$x_{ab}$ with $a\neq b$ of the following form:
 
$$ 
\begin{array}{llll}  x_{ik}x_{hj}  & \mbox{if \ \  } & i=j \mbox{\ \ or \
    \  } h=k & (1)  \\ 
 x_{ij}x_{hk} &  \mbox{if \ \  } & i<h \mbox{\ \ and \ \  } j<k     & (2) 
\end{array} 
\eqno{(*)}$$ Note that (1) is obvious by construction while (2)
follows immediately from the fact that if $i<h$ and $j<k$ then
$\max(|i-j|, |h-k|)<\max(|k-i|, |h-j|)$.

\begin{proposition}
\label{inigeni} The ideal $H(2,n)$ generated by the monomials described in 
$(*)$ is an initial ideal of $I(2,n)$ with respect to $\succ$.
\end{proposition}

It is clear that $H(2,n)\subseteq \ini_\succ(I(2,n))$. To prove
equality we use the following well-known fact.

\begin{lemma}\label{kntrick}  Let $J$ and $I$ be homogeneous ideals in a
  polynomial ring $R$. Assume that $J\subseteq I$, $\dim R/J=\dim
  R/I$, $\deg R/J \geq \deg R/I$ and $J$ is pure (i.e. all its
  associated primes have the same dimension). Then $J=I$.
\end{lemma} 

The proof of this fact is a simple exercise in primary decompositions.
Suppose $d=\dim R/I=\dim R/J$. Let $J=Q_1\cap \dots \cap Q_s$ be the
primary decomposition of $J$. By assumption $\dim R/Q_i=d$ for all
$i$.  Then $\deg R/J=\sum \deg R/Q_i$. Now, since $J\subseteq I$, each
primary component of $I$ of dimension $d$ must contain one of the
$Q_i$. As $\deg R/I=\deg R/J$, this forces the intersection of the
primary components of $I$ of dimension $d$ to be exactly $J$. So
$I\subseteq J$ and hence $I=J$.

We apply Lemma~\ref{kntrick} with $I=\ini_\succ(I(2,n))$ and
$J=H(2,n)$.  Because passing to initial ideals is a flat deformation
the dimension and the degree of $\ini_\succ(I(2,n))$ are
equal to that of $I(2,n)$: 
$\dim K[x_{ij}]/I(2,n) =2n-1$ and $\deg K[x_{ij}]/I(2,n)= { 2n-2
  \choose n-1}$.  
So to prove Proposition \ref{inigeni}  it suffices to show the following. 

\begin{lemma} \label{puredd}  
The ideal  $H(2,n)$ is pure of dimension  $2n-1$  and 
degree ${ 2n-2 \choose n-1}$.
\end{lemma}

\begin{proof} Let  $\Delta$ be the simplicial complex associated with
  $H(2,n)$. By construction the cone points of $\Delta$ are $CP=\{
x_{11},\dots,x_{nn} \}$ and we may concentrate our attention on
$\Delta' := \core(\Delta)$.  We have to show that $\Delta'$ is pure
and has exactly ${ 2n-2 \choose n-1}$ facets of dimension $n-2$. Let
us describe the facets of $\Delta'$.  For every nonempty proper subset
$R$ of $[n]$ we define:
$$\Delta_R=\{ F \in \Delta' : F\subseteq  R\times ([n]\setminus R)\}.$$
The generators of $H(2,n)$ of type (1)
imply that every face $F=\{ (a_1,b_1), \dots, (a_k,b_k)\}$
of $\Delta'$ has $\{a_1,\dots, a_k\}\cap \{b_1,\dots,
b_k\}=\emptyset$. In particular $F$ belongs to $\Delta_R$ with
$R=[n]\setminus \{b_1,\dots, b_k\}$, and hence $\Delta'= \cup
\Delta_R$. The generators of type (2) imply that $\Delta_R$ is exactly
the simplicial complex of the subsets of the grid $R\times
([n]\setminus R)$ which do not contain $2$-diagonals. 
If $R=\{r_1,\dots, r_p\}$ and
$[n]\setminus R=\{c_1,\dots, c_{n-p}\}$ with $r_1<\dots<r_p$ and
$c_1<\dots<c_{n-p}$, then a facet of $\Delta_R$ is a path in the grid
$R\times ([n]\setminus R)$ from $(r_p,c_1)$ to 
$(r_1,c_{n-p})$. We deduce two important facts. First, any facet of
$\Delta_R$ has $n-1$ elements and it involves all the elements of $R$
as row indices and all the elements in $[n]\setminus R$ as column
indices. Second, a facet of $\Delta_R$ cannot be a
facet of $\Delta_{S}$ if $R\neq S$. So the set of facets of $\Delta'$
is simply the disjoint union of the facets of $\Delta_R$ as
$R$ varies.  This implies that $\Delta'$ is pure of dimension
$n-2$. Note that the number of facets of $\Delta_R$ is ${ n-2 \choose
p-1}$ ($p = |R|$). In general, the number of paths in a
$a \times b$ grid from the bottom left to the top right is ${a+b-2
\choose a-1}$. Then the number of facets of $\Delta'$ is:
$$\sum_{p=1}^{n-1} {n \choose p}{n-2 \choose p-1} = \sum_{p=0}^{n-2}
{n \choose n-1-p}{n-2 \choose p} = { 2n-2 \choose n-1}.$$
\end{proof} 

In order to prove that $H(2,n)$ is Gorenstein, according to Lemma
\ref{lem1}, it suffices to prove that every face of $\Delta'$ of
codimension one is contained in exactly two facets and we need to
describe a shelling.  Actually we will describe a two-way shelling of
$\Delta'$.  First some notation.

Given a grid of size $a\times b$ we look at paths connecting the lower
left corner box $S$ (start) to the upper right corner box $E$ (end)
consisting of horizontal steps to the right or vertical steps up. Such
a path consists of $4$ types of points as we go from $S$ to $E$: a
left turn ($\lt$), a right turn ($\rt$), isolated point in a column
($\bullet$), and isolated point in a row ($\circ$).  This definition
is illustrated by Table \ref{typoints}.

\begin{scriptsize}
\begin{table}[ht]
\caption{}
\renewcommand\arraystretch{1.5}
\noindent
\[ \ \  \mbox{ A path } \quad
\begin{array}{|c|c|c|c|c|}
\hline    &    &     & *  & * \\
\hline    &  * &  *  & *  &  \\
\hline    &  * &    &     &  \\
\hline  * &  * &    &     &  \\
\hline
\end{array}
 \quad  \mbox{ and the type of its points }  \quad
\begin{array}{|c|c|c|c|c|}
\hline    &     &     &  \rt  & \bullet\\
\hline    &  \rt &  \bullet  & \lt &  \\
\hline    &  \circ &     &     &  \\
\hline  \bullet & \lt &     &     &  \\
\hline
\end{array}
\] 
\label{typoints}
\end{table} 
\end{scriptsize}
We say that a subset $A$ of the points of the grid has full
support if it intersects each row and each column.  Clearly a path
from $S$ to $E$ has full support.

\begin{lemma}\label{x3} Let  $P$ be a path in a grid  and $x\in P$.  We have: 
\begin{itemize}  
\item[i)] $x$ is a turn ($\lt$ or $\rt$) of $P$ if and only if
  $P\setminus \{x\}$ has full support if and only if there is exactly
  one other path $Q$ in the grid containing $P\setminus \{x\}$. The path
  $Q$ is obtained from $P$ and $x$ by ``flipping" $x$.
\item[ii)] $x$ is of type $\bullet$  or $\circ$ in $P$ if and 
  only if $P\setminus \{x\}$ does not have full support 
if and only if $P$ is the only
path in the grid containing $P\setminus \{x\}$.
\end{itemize}
\end{lemma}

To give an example, flipping the turn on the last row and second
column in the path of Table \ref{typoints} we get the path of Table
\ref{flip}.
\begin{scriptsize}
\begin{table}[ht]
\caption{}
\renewcommand\arraystretch{1.5}
\noindent 
\[ \ \  
\begin{array}{|c|c|c|c|c|}
\hline    &    &     & *  & * \\
\hline    &  * &  *  & *  &  \\
\hline  * &  * &     &     &  \\
\hline  * &    &     &     &  \\
\hline
\end{array}
 \]
\label{flip}
\end{table} 
\end{scriptsize}
 
\begin{lemma}\label{x4}  Let $P \in \Delta_R$ be a facet of $\Delta'$
and let $x$ be a point of $P$.  Then there are exactly two facets $P$
  and $Q$ of $\Delta'$ containing $P\setminus \{x\}$. The path $Q$ is
  described as follows:
\begin{itemize}
\item[i)] If $x$ is a turn of $P$ then $Q$ is the path of the grid $R
  \times ([n]\setminus R)$ (i.e. a facet of $\Delta_R$) obtained by
  flipping $x$.
\item[ii)] If $x$ is of type $\bullet$ then let $c$ be the column index of
  $x$. Set $R'=R\cup \{c\}$. Then $P\setminus \{x\}$ is a face of
  $\Delta_{R'}$ contained in a unique facet $Q$ of $\Delta_{R'}$. (See
  Table~\ref{thisQ}).
\item[iii)] If $x$ is of type $\circ$ then let $r$ be the row index of $x$.
  Set $R'=R\setminus \{r\}$. Then $P\setminus \{x\}$ is a face of
  $\Delta_{R'}$ contained in a unique facet $Q$ of $\Delta_{R'}$.
\end{itemize}
\end{lemma}

\begin{proof} i)  Note that  $P\setminus \{x\}$ has full support in
  the grid $R\times ([n]\setminus R)$. Hence  any facet of 
  $\Delta'$ containing $P\setminus \{x\}$ is a facet of $\Delta_R$.
  Now we use  i) of Lemma \ref{x3}. \\
ii) The support of $P\setminus \{x\}$ is $R\times ( [n]\setminus
R')$. So, among all the $\Delta_S$, $P\setminus \{x\}$ belongs only to
$\Delta_R$ and to $\Delta_{R'}$.  In both $\Delta_R$ and $\Delta_{R'}$
the set $P\setminus \{x\}$ does not have full support.  By ii) of Lemma 
\ref{x3} there is exactly one facet $P$ in $\Delta_R$ and exactly one facet 
$Q$ in$\Delta_{R'}$ containing $P\setminus \{x\}$. 
The statement iii) is dual to statement ii).
\end{proof}

For an illustration of the construction of Lemma~\ref{x4} ii) see
Table \ref{thisQ}, where $n=9$, $R=\{1,4,6,9\}$, $c=5$, $x=(4,5)$. The first
two arrays show $P$ and $P\setminus \{x\}$ in the grid
$R\times([n]\setminus R)$ and the second two show $P\setminus \{x\}$
and $Q$ in $R'\times([n]\setminus R')$.

\begin{scriptsize}
\begin{table}[ht]
\caption{}
\renewcommand\arraystretch{1.5}
\noindent 
\[  \begin{array}{rlrl}  P= & 
\begin{array}{|c|c|c|c|c|c|}
\hline    &  2 &  3 &  5  & 7  & 8 \\ 
\hline 1  &    &    &     & *  & * \\
\hline 4  &    &  * &  \bullet  & *  &   \\
\hline 6  &  * &  * &     &    &    \\
\hline 9  &  * &    &     &    &   \\
\hline
\end{array}
 \quad  \to \quad &  P\setminus\{x\}= &
\begin{array}{|c|c|c|c|c|c|}
\hline    &  2 &  3 &  5  & 7  & 8 \\ 
\hline 1  &    &    &     & *  & * \\
\hline 4  &    &  * &     & *  &   \\
\hline 6  &  * &  * &     &    &    \\
\hline 9  &  * &    &     &    &   \\
\hline
\end{array} 
 \quad \to \\ \\  P\setminus\{x\} =&
 \begin{array}{|c|c|c|c|c|}
\hline    &  2 &  3 & 7  & 8 \\ 
\hline 1  &    &    & *  & * \\
\hline 4  &    &  * & *  &   \\
\hline 5  &    &    &    &   \\
\hline 6  &  * &  * &    &    \\
\hline 9  &  * &    &    &   \\
\hline
\end{array}
 \quad \to \quad  &   Q=& 
\begin{array}{|c|c|c|c|c|}
\hline    &  2 &  3 & 7  & 8 \\ 
\hline 1  &    &    & *  & * \\
\hline 4  &    &  * & *  &   \\
\hline 5  &    &  \circ   &    &   \\
\hline 6  &  * &  * &    &    \\
\hline 9  &  * &    &    &   \\
\hline
\end{array}
\end{array}
 \]
\label{thisQ}
\end{table} 
\end{scriptsize} 

Now we describe the shelling. First we order the set of
nonempty  proper subsets of $[n]$. Such a subset is
represented as a strictly increasing sequence of integers. 

$$S= \{a_1,\dots, a_s \} < R= \{b_1,\dots, b_t \} \iff \left\{
\begin{array}{ll}  a_j<b_j  \mbox{ for the smallest } j \mbox{ such
    that } \\
a_j\neq b_j \\ 
\mbox{ or }   \\ s<t \mbox { and  } a_i=b_i \mbox{ for all } i=1,\dots, s 
\end{array} \right.
$$

\begin{definition}\label{she1} Let  $F$ and $G$ be facets of
  $\Delta'$, say $F$ is a facet of $\Delta_R$ and $G$ is a facet of
  $\Delta_S$.  We set:

$$F<G \iff  \left\{ 
\begin{array}{ll}  R<S  \\
\mbox{ or }   \\ R=S  \mbox { and  }  F<G \mbox{ in the standard
  shelling of  } \Delta_R 
\end{array} \right.   
$$
\end{definition} 

The standard shelling of $\Delta_R$ is defined as follows: let $F,G$
be facets (paths) in the corresponding grid. Then we set $F<G$ if the
first step in which they differ (always going from bottom-left to
top-right) is vertical for $F$ and (hence) horizontal for $G$. See
Table \ref{ts33} for the standard shelling in the $3\times 3$ grid.

\begin{scriptsize}
\begin{table}[ht]
\caption{}
\renewcommand\arraystretch{1.5}
\[ 
\begin{array}{l} 
\begin{array}{|c|c|c|}
\hline  *  &  * &   *   \\
\hline  *  &    &       \\
\hline  *  &    &       \\
\hline
\end{array}
\, \, < \,\,
\begin{array}{|c|c|c|}
\hline     &  * &   *   \\
\hline  *  &  * &       \\
\hline  *  &    &       \\
\hline
\end{array}
\,\, < \,\, 
\begin{array}{|c|c|c|}
\hline     &    &   *   \\
\hline  *  &  * &   *    \\
\hline  *  &    &       \\
\hline
\end{array} 
\,\, < \,\, 
\begin{array}{|c|c|c|}
\hline     &  * &   *   \\
\hline     &  * &        \\
\hline  *  &  * &       \\
\hline
\end{array}
\,\, < \,\,
\begin{array}{|c|c|c|}
\hline     &    &   *   \\
\hline     &  * &   *     \\
\hline  *  &  * &       \\
\hline
\end{array}
\,\, < \,\, 
\begin{array}{|c|c|c|}
\hline     &    &   *   \\
\hline     &    &   *   \\
\hline  *  &  * &   *   \\
\hline
\end{array}
\end{array}
\]
\label{ts33}
\end{table}  
\end{scriptsize}

For every facet $F$ of $\Delta_R$ we define: 

$$F^-=\begin{array}{l} 
\{ x\in F : x \mbox{ is a left  turn    } \} \,\cup \\
 \{ x\in F : x \mbox{ is of type $\bullet$ and its column index is } <\max(R)
 \} \,\cup \\ 
 \{ x \in F : x \mbox{ is of type $\circ$ and its row index is} \max(R) \}
\end{array}
$$ 

and 

$$F^+=F\setminus F^-=
\begin{array}{l} 
\{ x\in F : x \mbox{ is a right turn } \} \, \cup \\
 \{ x\in F : x \mbox{ is of type $\bullet$  and its column index is }
 >\max(R) \} \, \cup \\ 
 \{ x \in F : x \mbox{ is of type $\circ$ and its row index is}<\max(R) \}
\end{array}
$$ 

In  Table \ref{+-} the symbols  + or --  mark whether that point is in
$F^+$ or $F^-$.  
\begin{scriptsize}
\begin{table}[ht]
\caption{}
\renewcommand\arraystretch{1.5}
\noindent 
\[   
\begin{array}{|c|c|c|c|c|c|}
\hline    &  2 &  3 &  5  & 7  & 10 \\ 
\hline 1  &    &    &     & +  & + \\
\hline 4  &    &  +&  -& -  &   \\
\hline 6  &  +&  -&     &    &    \\
\hline 8  &  +&    &     &    &   \\
\hline 9  &  -&    &     &    &   \\
\hline
\end{array}
\]
\label{+-} 
\end{table} 
\end{scriptsize}

\begin{proposition}\label{xpro} The total order of the facets of
  $\Delta'$ in Definition~\ref{she1} is a two-way shelling of
$\Delta'$. Precisely, for every facet $F$ of $\Delta'$ one has:
$$\langle F  \rangle \cap \langle G : G<F  \rangle =\langle F\setminus
\{x\} : x\in F^-  \rangle \eqno{(1)}$$  and  
$$\langle F \rangle \cap \langle G : G>F \rangle =\langle F\setminus
\{x\} : x\in F^+ \rangle. \eqno{(2)} $$ where in (1) it is assumed
that $F$ is not the minimal facet of $\Delta'$ and in (2) it is
not the maximal.
\medskip
\end{proposition}

In order to prove Proposition \ref{xpro} we will show the two 
inclusions $\supseteq$ and $\subseteq$ in (1) and (2) separately.
The first inclusion is equivalent to the following

\begin{claim}\label{cl1}   For every facet $F$ of $\Delta'$ and for
  $x\in F$ let $G$ be the unique facet other than $F$ containing
$F\setminus\{x\}$.  Then we have $G>F$ if $x\in F^+$
and $G<F$ if $x\in F^-$. \end{claim}
\begin{proof}Suppose  $F$ is a facet of
$\Delta_R$. The statement is
clear if $x$ is a turn where $G$ is obtained by flipping the turn $x$. 
In this case,  $G<F$ if $x$ is a left turn and $G>F$ if $x$ is a right turn.  
If $x$ is of type $\bullet$ then $G\in \Delta_{R'}$ where $R'=R\cup \{c\}$ and
$c$ is the column index of $x$. We conclude that $R'<R$ if and only if
$c<\max(R)$. Finally if $x$ is of type $\circ$ then $G\in \Delta_{R'}$
where $R'=R\setminus \{r\}$ and $r$ is the row of $x$. 
And this time we conclude that $R'<R$ if and only if $r=\max(R)$. 
\end{proof}

The reverse inclusions in (1) and (2) translate to two more claims. 

\begin{claim}\label{cl2}  If $F, H$ are facets of $\Delta'$  and $H<F$
  then there exists $x\in F^-$ such that $x\not\in H$.
\end{claim}

\begin{claim}\label{cl3}  If $F,G$ are facets   of $\Delta'$  and
  $F<G$  then there exists $y\in F^+$ such that $y\not\in G$.
\end{claim}

\begin{proof}[Proof of Claim~\ref{cl2}] Suppose $F$ is a facet of $\Delta_R$.
  If $H$ also belongs to $\Delta_R$, then the desired $x$ is a left
turn of $F$. If instead $H$ is in $\Delta_S$ for some $S\neq R$ then
$S<R$ because $H<F$. We let  $S=\{a_1 < \dots <a_s\}$ and
$R=\{b_1 < \dots <b_t\}$. There are two cases.  \smallskip
 
\noindent Case 1: $a_j<b_j$ for some $j$ and $a_i=b_i$ for every
$i<j$. Then $a_j\not\in R$ and $a_j<b_j\leq \max(R)$. 
The points of $F$ in column $a_j$ are
not in $H$ since $a_j$ is a row index for $H$. So it is enough to show
that column $a_j$ intersects $F^-$.  If $F$ has an isolated point $x$ in
this column, then we are done since we know that $a_j<\max(R)$.
If $F$ has a left turn in column $a_j$ then we are also done.  
Otherwise $a_j$ is the first column of the
grid of $F$ and the first step of the path is vertical. But the
starting point of the path is of type $\circ$ in the row with index
$\max(R)$ and column $a_j$.  This concludes the proof in this case.
\smallskip
 
\noindent Case 2: $s<t$ and $a_i=b_i$ for $i=1,\dots,s$. 
The points of $F$ in row $b_t$ are not in $H$ since
$b_t$ is a column index for $H$. So it is enough to show that row
$b_t$ of $F$ intersects $F^-$. This is clear because in the last row
we have either a left turn or an element of type $\circ$ (note that $F$ has
at least two rows).
\end{proof}
 
\begin{proof}[Proof  of Claim~\ref{cl3}] As above if $G$ is a facet 
of $\Delta_R$ such that $F \in \Delta_R$ then the desired $y$ is 
a right turn of $F$. If instead $G$ is in $\Delta_S$ for
some $S\neq R$ then $S>R$ because $G>F$. We let  $S= \{a_1 < \dots <a_s\}$
and $R=\{b_1 < \dots < b_t\}$ and study two cases.
\smallskip
 
\noindent Case 1: $a_j>b_j$ for some $j$ and $a_i=b_i$ for every
$i<j$. The points of $F$ in row $b_j$ are not in $G$ since 
$b_j$ is a column index for $G$.  So we
are done if row $b_j$ intersects $F^+$.  This is the case if $F$ has a
right turn in row $b_j$.  It is also the case if $F$ has an isolated
point in row $b_j$ and $b_j<\max(R)$.  So we may assume that
$b_j=\max(R)$ and $F$ has no right turn in that row.  Now 
either $j=1$ (i.e.  $R=\{b_1\}$) or $j>1$ and
the first step of $F$ is vertical.  If $j=1$ then $(1,a_1)$
is of type $\bullet$ for $F$ and we are done. If $j>1$ and the first step of
$F$ in vertical then $a_j$ is a column index for $F$ and is
$>\max(R)=b_j$.  If in column $a_j$ for $F$ we have an isolated point
or a right turn then we are done. There is just one
possibility left: $a_j$ is the last column index for $F$ and there is no
point of type $\bullet$ in that column.  So the ending point of $F$ is
reached with a vertical step. But then the ending point $(b_1,a_j)$ is
of type $\circ$ and $b_1<b_j$. So we are done. \smallskip

\noindent Case 2: $t<s$ and $a_i=b_i$ for $i=1,\dots,t$.
All the points of $F$ in column $a_s$
are not in $G$ since $a_s$ is a row index for $G$. So we are done if
column $a_s$ intersects $F^+$. This is the case if $F$ has a point of
type $\bullet$  in column $a_s$ or a right turn in that column. Otherwise
$a_s$ must be the largest column index for $F$ and the last step of
$F$ is vertical. But then $t>1$ (because there is a vertical step) and
the last point of $F$, namely $(b_1,a_s)$ is of type $\circ$ for $F$ with
row index $b_1<\max(R)=b_t$. So $(b_1,a_s)$ is in $F^+$.  This
concludes the proof of Claim \ref{cl3}.
\end{proof} 
 
Thus we have shown that $H(2,n)$ gives a Gorenstein initial complex of
$I(2,n)$.  The goal of the rest of this section is to prove that for
many reverse lexicographic initial ideals similar to $H(2,n)$ the core
of the initial complex is the boundary of a simplicial polytope. To
this end we construct a reverse lexicographic triangulation of the
point configuration ${\mathcal P}$ whose toric ideal is $I_{\mathcal
P}=I(2,n)$. We first review some facts about these triangulations and
${\mathcal P}$.

Let $\mathcal{P} = \{a_1, \ldots, a_n\} \subset \ZZ^d$ be a point 
configuration. The reverse lexicographic triangulation of ${\mathcal P}$ 
(as well as the corresponding polytope $P$) with respect to 
the ordering $a_1 \succ a_2 \succ \cdots \succ a_n$ is obtained as follows 
(see \cite[Chapter 8]{Stu1}):  let $F_1, \ldots, F_k$ be the facets of 
$P$ that do not contain $a_n$. Then 
$$\Delta_\succ(\mathcal{A}) \,\, = \,\, \bigcup_{i=1}^k \bigcup_{G \in
\Delta_\succ(F_i)} G \cup \{a_n\}$$ where $G$ runs over the facets of
$\Delta_\succ(F_i)$. Observe that the definition implies that $a_n$ is
a cone point of $\Delta_\succ(P)$.

Denote by $I(2,m,n)$ the toric ideal of $\Segre(m,n)$ generated by the
$2$-minors of a generic $m \times n$ matrix. In this case, ${\mathcal
P}$ is $\Sigma_{m-1} \times \Sigma_{n-1}$ where $\Sigma_k$ is the
standard simplex in $\RR^{k+1}$ of dimension $k$. In this case the
point configuration is
$$ \mathcal{P}(m,n) \,\, := \,\, \{e_i \oplus f_j \, : \, 1 \leq i
\leq m, 1 \leq j \leq n\}$$ where $e_i$ and $f_j$ are the standard
unit vectors in $\RR^m$ and $\RR^n$, respectively. We will identify
the columns of $\mathcal{P}(m,n)$ with the variables in the polynomial
ring $K[x_{ij}]$.  The convex hull of $\mathcal{P}(m,n)$ which we
denote by $P(m,n)$ has dimension $m+n-1$. Each face of $P(m,n)$ is $F
\times G$ where $F$ and $G$ are faces of $\Sigma_{m-1}$ and
$\Sigma_{n-1}$, respectively. In other words, facets of $P(m,n)$ are
$F \times \Sigma_{n-1}$ and $\Sigma_{m-1} \times G$ where $F$ and $G$
run over the facets of $\Sigma_{m-1}$ and $\Sigma_{n-1}$,
respectively.  This implies the following.

\begin{proposition} Let $u_i$ and $v_j$ be the coordinate functions
of $\RR^m$ and $\RR^n$ ($m , n \geq 2$). Then the facets of $P(m,n)$
in $\RR^m  \oplus \RR^n$ are precisely the $m+n$ faces supported by
$u_i = 0$ for $i=1, \ldots, m$ and $v_j = 0 $ for $j=1, \ldots,n$.  
\end{proposition}

We will also need  the following lemma.

\begin{lemma} \label{corelemma} Assume that $0 \leq i \leq m \leq n$. 
Let $x_{11} \prec x_{22} \prec \cdots \prec x_{ii} \prec \{ x_{ij}\, : \, i
\neq j \}$ be a reverse lexicographic order where the variables in the
latter set are ordered arbitrarily.  Then $x_{11}, \ldots, x_{ii}$ are
cone points of the triangulation $\Delta_\succ$. Moreover, if $i=m >
1$, the simplicial complex obtained by removing $x_{11}, \ldots,
x_{mm}$ is $\core \Delta_\succ$.
\end{lemma}

\begin{proof} We induct on $m+n$. The first non-trivial
cases are $m+n = 3$ and $m+n=4$, and the statements are easy to check.
The case $i=0$ is vacuous for any $m+n$. So suppose $i \geq 1$. By the
definition of $\Delta_\succ$ we know $x_{11}$ is a cone point. There
are exactly two facets of $P(m,n)$ that do not contain $x_{11}$,
namely the facets defined by $u_1 =0$ and $v_1 = 0$. These facets are
isomorphic to $P(m-1,n)$ and $P(m,n-1)$ respectively. They go with
$I(2,m-1,n)$ and $I(2, m,n-1)$ corresponding to generic matrices
obtained by deleting the first row and deleting the first column
(respectively) of an $m\times n$ matrix.  In the first case,
by cyclically permuting the columns, and in the second case by
cyclically permuting the rows, we will be in the case $m+n-1 < m +n $
and $x_{22} \prec \cdots \prec x_{ii}$ are the variables that are
smallest. By induction they are cone points on both facets, and hence
cone points of $\Delta_\succ$. For the last statement, observe that
after removing $x_{11}, \ldots, x_{mm}$ the remaining faces that need
to be triangulated are defined by $u_i = 0$ for $i \in I \subset [m]$
together with $v_j = 0$ for $j \in [m] \setminus I$.  If there were
another cone point $x_{ij}$, the corresponding $e_i \oplus f_j$ had to
be in every one of these faces. But clearly that cannot happen.
\end{proof}

\begin{theorem} \label{spherethm}
Assume that $0 \leq i \leq m \leq n$.  Let $x_{11} \prec x_{22} \prec
\cdots \prec x_{ii} \prec \{ x_{ij}\, : i \neq j \}$ be a reverse
lexicographic order.  After removing the cone points $x_{11}, \ldots,
x_{ii}$ from the triangulation $\Delta_\succ$, the remaining
simplicial complex is a $m+n-i-2$ dimensional ball if $m<n$ or if
$m=n$ and $i < m$, and it is a $n-2$ dimensional sphere if $i=m=n$.
\end{theorem} 

\begin{proof} Again, the proof is by induction on $m+n$. One more time
the cases $m+n=3$ and $m+n=4$ are easy to check. As in the proof of
the above lemma, after removing the cone point $x_{11}$, the rest of
the triangulation is the union of the reverse lex triangulations of
the two facets of $P(m,n)$ defined by $u_1= 0$ and $v_1=0$
respectively. These facets were isomorphic to $P(m-1,n)$ and
$P(m,n-1)$, and we we will use our induction hypothesis on
them. Assume that we remove the cone points $x_{22}, \ldots, x_{ii}$
from these two facets. We get the following statements by induction:
if $m<n$ or if $i<m=n$, the two simplicial complexes are $m+n-i-2$
dimensional balls.  They are glued along the simplicial complex
obtained by triangulating the unique face at the intersection of the
two facets, namely the face defined by $u_1=v_1=0$, and removing the
cone points $x_{22},\ldots, x_{ii}$. This face is isomorphic to
$P(m-1,n-1)$, and hence after the removal we get a
$m+n-i-3$-dimensional ball. But this one-lower-dimensional ball is on
the boundary of the two balls. So the gluing gives again an
$m+n-i-2$-dimensional ball. When $i=m=n$, we obtain two
$n-2$-dimensional balls, glued by an $n-3$-dimensional sphere. If we
can show that this $n-3$-sphere is exactly the boundary of the two
balls, then after gluing we will get a $n-2$-dimensional sphere. Let's
concentrate on one ball $B$, obtained from the facet $u_1=0$. After
removing the cone points, this simplicial complex is the union of
simplicial complexes obtained by triangulating the faces $F_I$ defined
by $u_1 =0$ and $u_i =0$ for $i \in I \subset \{2,\ldots, m\}$ and
$v_j =0 \in J = \{2,\ldots, m\} \setminus I$. So the simplices that
will make up the boundary of $B$ are precisely the simplices on the
facets of $F_I$ which belong to a unique $F_I$. The facets of $F_I$
are obtained by either setting $u_s = 0$ where $s \in J$ or setting
$v_t = 0$ where $t \in I \cup \{1\}$. In the first case, this facet of
$F_I$ is also a facet of $F_{I \cup \{s\}}$ defined by $v_s = 0$. In
the second case, if $t \neq 1$, it is the facet of $F_{I \setminus
\{t\}}$ defined by $u_t = 0$. Only when $t=1$, this facet of $F_I$
belongs to the $n-3$-dimensional sphere which is on the boundary of
this $B$. Symmetric arguments hold for the second facet defined by
$v_1 = 0$, and we are done.
\end{proof}

\begin{figure} 
\epsfig{file=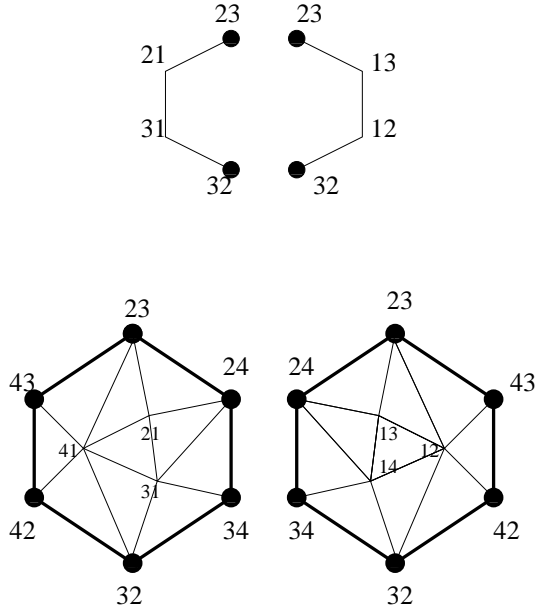, height={8cm}}
\caption{How the balls glue to a  sphere}
\label{glue}
\end{figure}

\begin{corollary} Any initial ideal of $I(2,n)$ with
respect to the reverse lex order $x_{11} \prec x_{22} \cdots \prec
x_{nn} \prec \{ x_{ij} \, : \, i \neq j \}$ is Gorenstein and
quadratic.
\end{corollary} 
\begin{proof} Any initial ideal of $I(2,m,n)$  is
squarefree, since ${\mathcal P}(m,n)$ is a totally unimodular 
configuration. So it is enough to show that $\core \Delta_\succ$ is 
a simplicial sphere. And this follows from Lemma \ref{corelemma}
and Theorem \ref{spherethm}. The fact that any such initial ideal is  
quadratic is a consequence of the following remark. 
\end{proof}

\begin{remark} 
A set of polynomials is a universal Gr\"obner basis for the ideal they
generate if it is a Gr\"obner basis of the ideal with respect to all
term orders, and it is a reverse lex universal Gr\"obner baisis if it
is a Gr\"obner basis with respect to all reverse lex term orders.  The
$2$-minors of a generic $m\times n$ matrix do not form a universal
Gr\"obner unless $\min(m,n)=2$. A universal Gr\"obner basis for the
ideal $I(2,m,n)$ can be described in terms of the cycles of the
complete bipartite graph $K_{m,n}$, see \cite[4.11 and 8.11]{Stu1} or
\cite[8.1.10 ]{Vil}.  Nevertheless, the $2$-minors of a generic
$m\times n$ matrix do form a universal reverse lex Gr\"obner basis of
$I(2,m,n)$.  This can be checked by using the Buchberger criterion.
\end{remark}

Let's illustrate the theorem for $m=n=2$, $m=n=3$ and $m=n=4$.  In the
first case $\core \Delta_\succ$ consists of the two isolated vertices
$x_{12}$ and $x_{21}$, and this is a $0$-sphere. For the other cases
Figure~\ref{glue} shows the two balls and how they glue along their
boundary. In the last case, we order the variables so that the initial
terms of the minors that do not touch the variables $x_{11}, \ldots,
x_{44}$ are the main diagonals. 
 

\begin{corollary} The simplicial sphere constructed
in Theorem \ref{spherethm} is the boundary of a simplicial polytope. 
\end{corollary}
\begin{proof}
The simplex spanned by $e_i \oplus f_i$ corresponding to $x_{11},
\cdots, x_{nn}$ is a special simplex as in Section \ref{gortoric}, and
the reverse lexicographic term order we used is the kind in Lemma
\ref{lemma3.4Athan}.  Now the result follows from this lemma and
Theorem \ref{thm3.5Athan}.
\end{proof}

\section{Gorenstein Veronese varieties}\label{Veronese} 

In Section \ref{gortoric} we have indicated that the ideal defining
the Veronese variety $\Ver(r,n)$ is Gorenstein if and only if $r$
divides $n$. Theorem~\ref{BRcor7} guarantees that this ideal has a
squarefree Gorenstein initial ideal. In this section we look at the
case $r=2$, and give two independent proofs of the same result. These
results are in the same spirit as in Section \ref{Segre}: the first
constructs a squarefree initial ideal that corresponds to a simplicial
complex with a two-way shelling, and the second constructs an initial
complex which is a polytopal sphere.
   
Recall that the ideal $J(2,n)$ generated by the $2$-minors of an
$n\times n$ generic symmetric matrix $X=(x_{ij})$ is the defining
ideal of $\Ver(2,n)$.
It is an ideal of the polynomial ring $K[x_{ij}]=K[x_{ij} : 1\leq
i\leq j\leq n]$.  The Hilbert series of $K[x_{ij}]/J(2,n)$ is
$$\sum_i {n \choose 2i} z^i/(1-z)^{n}.$$
We assume that $n=2m$.  Then the degree of the $h$-vector is $m$ and the
$a$-invariant is $n-m = m $. Therefore any Gorenstein initial
complex of $J(2,n)$ must have exactly $m$ cone points.

The classical initial complex (associated to diagonal orders) is
described as follows: its facets are the paths in the ``upper
triangle" in an $n \times n$ grid
$$T=\{(i,j) \in \NN^2 :  1\leq i\leq j\leq n\}$$ with starting point
$(1,n)$ and end point $(i,i)$ for some $i$, $1\leq i\leq n$. 

Table \ref{thisisT} shows $T$ and a typical facet of the classical
initial complex of $J(2,4)$, where we put $\circ$ in those positions which
are in the triangle but not in the path.
\begin{scriptsize}
\begin{table}[ht]
\caption{}
\renewcommand\arraystretch{1.5}
\noindent
\[ 
\begin{array}{|c|c|c|c|c|}
\hline   \circ &  \circ  &  \circ   & \circ   \\
\hline     &  \circ  &  \circ   & \circ   \\
\hline     &     &  \circ   & \circ   \\
\hline     &     &      & \circ   \\
\hline
\end{array}  
\quad 
\begin{array}{|c|c|c|c|c|}
\hline   \circ &  \circ  &  \circ   & *   \\
\hline     &  \circ  &  *   & *   \\
\hline     &     &  *   & \circ   \\
\hline     &     &      & \circ   \\
\hline
\end{array} 
\] 
\label{thisisT}
\end{table}
\end{scriptsize}

The only cone point of the classical initial complex of $J(2,n)$ is
$(1,n)$, and hence it is not Gorenstein if $n>2$.  In order to
describe a Gorenstein initial complex we consider a term order such
that the initial term of a $2$-minor of $(x_ {ij})$ is its main
diagonal unless the main diagonal involves elements from the set
$CP=\{(1,n), (2,n-1), \dots, (m,m+1)\}$.  An example of such a term
order is a reverse lexicographic order where the variables
corresponding to $CP$ (namely $x_{1n}, x_{2,n-1}, \dots, x_{m,m+1}$)
are followed by the rest of the variables which are totally ordered so
that $x_{ij}\succ x_{hk}$ if $|i-j|<|h-k|$.  For instance, for $n=4$
this term order can be taken as the reverse lexicographic order with
$x_{11}\succ x_{22} \succ x_{33} \succ x_{44} \succ x_{12} \succ
x_{34} \succ x_{13} \succ x_{24} \succ x_{23} \succ x_{14}$.

By construction, the initial terms of the $2$-minors are the monomials
not involving variables in $CP$ of the following two kinds:
 
$$  
\begin{array}{ccc}  x_{ij}x_{hk}  & \mbox{  if  a+b=n+1 for some }
  a\in \{i,j\} \mbox{ and } b \in \{h,k\} & (1) \\ x_{ij}x_{hk} &
\mbox{ with } i\leq j, h\leq k, i<h, j<k & (2) \end{array}
\eqno{(**)}$$

Let $K(2,n)$ be the ideal generated by these monomials. We want to
show that $K(2,n)=\ini_\succ(J(2,n))$.  According to
Lemma~\ref{kntrick}, it suffices to show that $K(2,n)$ and
$\ini_\succ(J(2,n))$ have the same dimension and degree and that
$K(2,n)$ is pure. The dimension and the degree of $J(2,n)$ are $\dim
K[x_{ij}]/J(2,n) =n$ and $\deg K[x_{ij}]/J(2,n)= 2^{n-1}$.

\begin{lemma} \label{applykntrick}
Let $\Delta' = \core(\Delta)$ be the core of the simplicial complex
$\Delta$ associated with $K(2,n)$.  Then $\Delta'$ is pure and has
$2^{n-1}$ facets with $m$ vertices.
\end{lemma}
\begin{proof}
Consider the family $\AA$ of subsets $A$ of $[n]$ of cardinality $m$
and such that $i+j\neq n+1$ for every $i,j\in A$.  Note that any $A\in
\AA$ is completely determined by its intersection with $[m]$. In other
words, the cardinality of $\AA$ is $2^m$.
For any  $A\in \AA$ we set 
$$T_A=\{ (i,j) : i\leq j \mbox{ and } i,j \in A\} \qquad \mbox{and}
\qquad \Delta_A=\{ F \in \Delta' : F\subseteq T_A\}$$ The monomials of
type (1) imply $\Delta'=\cup_A \Delta_A$ as $A$ varies in $\AA$.  The
monomials of type (1) do not have any effect on $\Delta_A$ while those
of type (2) imply that $\Delta_A$ is exactly the simplicial complex of
the subsets of the small triangle $T_A$ which do not contain any
$2$-diagonal. In other words any $\Delta_A$ is the classical initial
complex of $J(2,m)$. Each $\Delta_A$ has $2^{m-1}$ facets each of
cardinality $m$. Each facet of $\Delta_A$ involves (either as a row or
column index) all the indices of $A$.  Then the set of the facets of
$\Delta$ is the disjoint union of the set of the facets of $\Delta_A$
with $A\in \AA$. It follows that $\Delta$ is pure and has $2^m2^{m-1}
=2^{n-1}$ facets.
\end{proof}

Our next goal is to prove that $\Delta$ is Gorenstein.  Given
$A=\{a_1,\dots,a_m\}$ with $a_1<\dots< a_m$, we consider paths in
$T_A$ starting with the box $S=(a_1,a_m)$ and ending with a box
$(a_i,a_i)$ on the diagonal.  Each step is either a horizontal step to
the left or a vertical step downwards. Such a path consists of $3$
types of points as we travel from $S$ to a diagonal box: a left turn
($\nlt$) with the convention that the last point is a left turn if the
last step to a diagonal box is horizontal; a right turn ($\nrt$) with
the convention that the last point is a right turn if the last step is
vertical; and an isolated point ($\bullet$) if this point is the only
one on the path on a row or column $a_j$. In the latter case we say
that the point is isolated with index $a_j$.  For an illustration of
this definition see Table \ref{stypes}.
\begin{scriptsize}
\begin{table}[ht]
\caption{} 
\renewcommand\arraystretch{1.5}
\noindent
\[ \ \  \mbox{ Path: } \qquad
\begin{array}{|c|c|c|c|c|c|}
\hline   0&  0&  0&  0& 0&  \phantom{r}*  \\
\hline    &  0&  0&  *& *&  \phantom{r}*  \\
\hline    &   &  0&  *& 0&  \phantom{r}0  \\
\hline    &   &   &  *& 0&  \phantom{r}0  \\
\hline    &   &   &   & 0&  \phantom{r}0  \\
\hline    \phantom{rt} &  \phantom{rt}  &  \phantom{rt}  &  \phantom{rt}  & 
\phantom{rt} &  \phantom{r}0  \\
\hline
\end{array}
 \qquad  \qquad \mbox{ Types:  } \qquad  
\begin{array}{|c|c|c|c|c|c|}
\hline   0&  0&  0&  0& 0&  \bullet  \\
\hline    & 0&  0&  \nlt&\bullet&  \nrt  \\
\hline    &   &  0&  \bullet& 0&  0  \\
\hline    &   &   &  \nrt& 0&  0  \\
\hline    &   &   &   & 0&  0  \\
\hline    &   &   &   &  &  0  \\
\hline
\end{array}
\]
\label{stypes}
\end{table} 
\end{scriptsize}
As in Lemma~\ref{x4} we prove:

\begin{lemma}\label{xyx}  Let $P\in \Delta_A$ be  a facet of $\Delta'$ and
 let $x$ be a point of $P$.  Then there are exactly two facets $P$ and
  $Q$ of $\Delta'$ containing $P\setminus \{x\}$.  The path $Q$ is
  described as follows:
\begin{itemize}
\item[i)] If $x$ is a turn of $P$ then $Q$ is the path of the triangle
  $T_A$ (i.e. a facet of $\Delta_A$) obtained by flipping $x$.
\item[ii)] If $x$ is of type $\bullet$ suppose 
 that it is the only point of the path
  involving the index $i\in A$. We  set 
$A'=A\setminus \{i\} \cup \{n+1-i\}$, and then
$P\setminus \{x\}$ is a face of
$\Delta_{A'}$ contained in a unique facet $Q$ of $\Delta_{A'}$. 
\end{itemize}
\end{lemma}   
 
We order the $A$'s lexicographically, i.e. if $A=\{a_1,\dots,a_m\}$ and
$B=\{b_1,\dots,b_m\}$ then

$$A< B \iff a_j<b_j \mbox{ for the smallest } j \mbox{ such that }
a_j\neq b_j.
$$
  
And also we define a total order on the set of
facets of $\Delta'$ which will turn out to be a shelling.

\begin{definition}\label{she2}  Let $F$ and $G$ be  facets of
  $\Delta'$,  say $F$ is a facet of $\Delta_A$ and $G$ is a facet of 
  $\Delta_B$.  We set:
$$F<G \iff  \left\{ 
\begin{array}{ll}  A<B  \\
\mbox{ or }   \\ A=B  \mbox { and  }  F<G \mbox{ in the standard
  shelling of  } \Delta_A 
\end{array} \right.   
$$
\end{definition} 

The standard shelling of $\Delta_A$ is defined as follows: let $F,G$
be facets (paths) in the corresponding $T_A$. Then $F<G$ if the first
step in which the paths differ going from top-right to bottom-left is
horizontal for $F$ and (hence) vertical for $G$. Table \ref{symsh44}
shows the standard shelling $\Delta_A$ when $m=4$.
\begin{scriptsize}
\begin{table}[ht]
\caption{}
\renewcommand\arraystretch{1.5}
\[  
\begin{array}{llll}
\begin{array}{|c|c|c|c|}
\hline  *  & *  &  *  & *  \\
\hline     &  \circ &  \circ  & \circ  \\
\hline     &    &  \circ  & \circ  \\
\hline     &    &     & \circ  \\
\hline
\end{array} &\quad < \quad
\begin{array}{|c|c|c|c|}
\hline  \circ  &  * &  *  & *  \\
\hline     &  * &  \circ  & \circ  \\
\hline     &    &  \circ  & \circ  \\
\hline     &    &     & \circ  \\
\hline
\end{array}  &\quad  < \quad
\begin{array}{|c|c|c|c|}
\hline  \circ  &  \circ &  *  & *  \\
\hline     &  * &  *  & \circ  \\
\hline     &    &  \circ  & \circ  \\
\hline     &    &     & \circ  \\
\hline
\end{array}  &\quad < \quad
\begin{array}{|c|c|c|c|}
\hline  \circ  &  \circ &  *  & *  \\
\hline     &  \circ &  *  & \circ  \\
\hline     &    &  *  & \circ  \\
\hline     &    &     & \circ  \\
\hline
\end{array} \quad < \\
\\
\begin{array}{|c|c|c|c|}
\hline  \circ  &  \circ &  \circ  & *  \\
\hline     &  * &  *  & *  \\
\hline     &    &  \circ  & \circ  \\
\hline     &    &     & \circ  \\
\hline
\end{array}  & \quad  < \quad
\begin{array}{|c|c|c|c|}
\hline  \circ  &  \circ &  \circ  & *  \\
\hline     &  \circ &  *  & *  \\
\hline     &    &  *  & \circ  \\
\hline     &    &     & \circ  \\
\hline
\end{array}  &\quad  < \quad
\begin{array}{|c|c|c|c|}
\hline  \circ  &  \circ &  \circ  & *  \\
\hline     &  \circ &  \circ  & *  \\
\hline     &    &  *  & *  \\
\hline     &    &     & \circ  \\
\hline
\end{array}  &\quad < \quad 
\begin{array}{|c|c|c|c|}
\hline  \circ  &  \circ &  \circ  & *  \\
\hline     &  \circ &  \circ  & *  \\
\hline     &    &  \circ  & *  \\
\hline     &    &     & *  \\
\hline
\end{array} 
\end{array}
\]\label{symsh44}
\end{table}  
\end{scriptsize}

\noindent
For every facet $F$ of $\Delta_A$ we define: 
$$F^-=
\{ x\in F : x \mbox{ is a right  turn    } \} \, \bigcup \,
\{ x\in F : x \mbox{ is  of type $\bullet$ with index } >m\} \mbox{ and}
$$ 
$$F^+=F\setminus F^-=
\{ x\in F : x \mbox{ is a left  turn    } \} \, \bigcup \,
\{ x\in F : x \mbox{ is  of type $\bullet$  with index } \leq m\}
$$ 
\noindent 
The proof of the next proposition is similar to that of
Proposition~\ref{xpro} but easier since here everything is fully
symmetric.
 
\begin{proposition}\label{ypro} The total order described above  is a
  two-way shelling of $\Delta'$. Precisely, for every facet $F$ of
  $\Delta'$ one has:

$$\langle F  \rangle \cap \langle G : G<F  \rangle =\langle F\setminus
\{x\} : x\in F^-  \rangle \eqno{(1)}$$  

and 
$$\langle F  \rangle \cap \langle G : G>F  \rangle =\langle F\setminus
\{x\} : x\in F^+  \rangle \eqno{(2)} $$ 
\end{proposition} 

By Lemmas~\ref{applykntrick}, \ref{xyx}, Proposition~\ref{ypro} and
applying Lemma \ref{lem1} we have proved that $K(2,n)$ is a squarefree
Gorenstein initial ideal of $J(2,n)$. A stronger statement would be to
claim that $\Delta'$ is the boundary complex of a simplicial polytope.
The rest of the section will prove this result.

The point configuration defining $\Ver(2,n)$ is
$${\mathcal P}(2,n) = \{e_i + e_j \, : 1 \leq i \leq j \leq  n\} \subset \RR^n,
$$ and the convex hull of these points is the polytope $P(2,n)$ with
vertices $2e_1, \ldots, 2e_n$. The variable $x_{ij}$ corresponds to
the point $e_i+e_j$. The facets of $P(2,n)$ are defined by the
coordinate hyperplanes $y_i = 0$ for $i=1,\ldots,n$.  To construct the
desired initial complex we will use the same reverse lexicographic
term order we introduced earlier.  We order the anti-diagonal
variables $AD = \{x_{m,m+1} \succ \cdots \succ x_{1n}\}$, and then we
order the rest of the variables $X^c = X \setminus AD$ in such a 
way, so that the faces of $P(2,n)$ not containing 
$e_1+e_n, \ldots, e_{m}+e_{m+1}$ are triangulated by unimodular simplices.
Because these faces are isomorphic to $P(2,m)$ such a coherent
ordering of $X^c$ are possible.

Now we construct the corresponding reverse lexicographic triangulation
of $P(2,n)$. Because $K(2,n)$ is a squarefree initial ideal, this
triangulation is unimodular. In fact we can describe its pieces
as we have done before: after ``pulling'' the points
corresponding to the anti-diagonal variables we are left with the
faces of $P(2,n)$ defined by setting $m$ coordinates $y_i = 0$ where  
$i \in A$ as in the proof of Lemma \ref{applykntrick}. These faces
are isomorphic to $P(2,m)$, and as Veronese polytopes they are
triangulated further using the usual diagonal term order into
unimodular simplices. We note two facts. First, after pulling
the point corresponding to $x_{1n}$ there are exactly two
facets $F_1$ and $F_n$ of $P(2,n)$, 
defined by $y_1=0$ and $y_n=0$ respectively,
which do not contain this point. Second, the intersection of
these two facets defined by $y_1 = y_n = 0$ is the polytope $P(2,n-2)$. 
\begin{theorem} \label{ver-sphere}
The core of the simplicial complex obtained by
triangulating $P(2,n)$ using the above term order
is an $(m-1)$-dimensional sphere.
\end{theorem}
\begin{proof} We use induction on $m$. The case for $m=1$ is clear since
$P(2,2)$ is the convex hull of $2e_1, e_1+e_2, 2e_2$ in $\RR^2$.
Suppose the statement is true for $k \leq m-1$.  As we observed above
the facet $F_1$ and $F_n$ of $P(2,2m)$ are triangulated polytopes, and
their intersection is $P(2,2m-2)$. Since these two facets contain the
remaining $m-1$ cone points in their intersection, the core of (the
triangulation of) $F_1 \cap F_2$ is the intersection of the cores of
$F_1$ and $F_n$.  By induction, the former is an $(m-2)$-dimensional
sphere.  The core of $F_1$ is supported on smaller Veronese polytopes
obtained by setting $m$ coordinates $y_i=0$ where $i \in A$ as in the
proof of Lemma \ref{applykntrick} and $n \not \in A$, and similarly
the core of $F_n$ is supported on those where $1 \not \in A$.  All of
the former contains the point $2e_n$ and all of the latter contains
the point $2e_1$.  Now since the core of $F_1 \cap F_2$ is supported
by faces obtained by setting both $y_1=y_n=0$, we conclude that the
core of the triangulation of $F_1$ is supported on cones with apex
$2e_n$ and that of $F_n$ is supported on cones with apex $2e_1$. Since
the core of $F_1 \cap F_n$ is an $m-2$-dimensional sphere we conclude
that the core of $F_1$ and $F_2$ are $(m-1)$-dimensional balls, and
their boundary is precisely the core of $F_1 \cap F_n$.  This shows
that the triangulation of $P(2,n)$ is an $(m-1)$-dimensional sphere.
\end{proof} 

\begin{figure} 
\epsfig{file=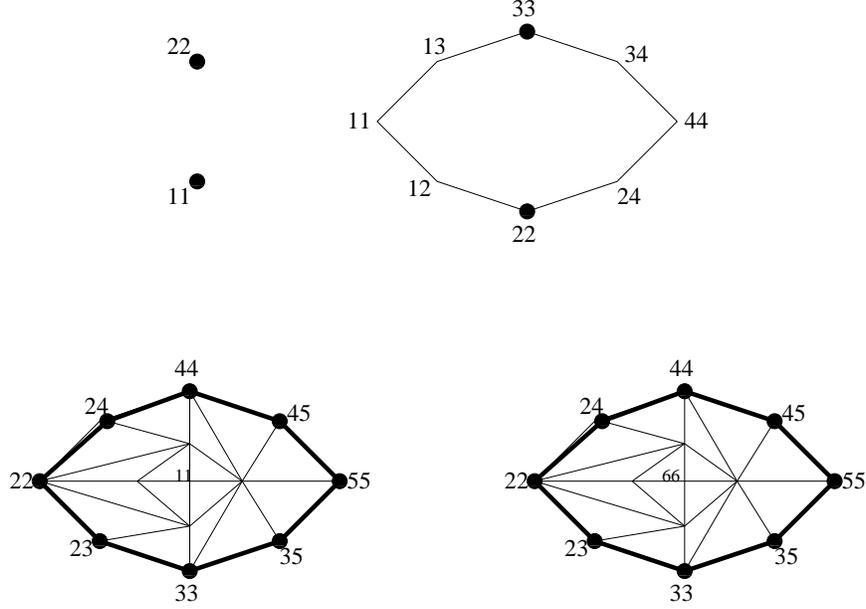, height={8cm}}
\caption{How the balls glue to a  sphere}
\label{glue2}
\end{figure}

Figure \ref{glue2} illustrates the construction in the above
proof for $m=1,2,3$.
\begin{corollary} The sphere constructed in Theorem \ref{ver-sphere}
is a polytopal sphere. 
\end{corollary}
\begin{proof} We note that the simplex that 
is the convex hull of the cone points $e_1+e_n, \cdots, e_m+e_{m+1}$
is a special simplex of $P(2,n)$ as defined in Section \ref{gortoric}.
The triangulation we get is a reverse lexicographic one 
used in Theorem \ref{thm3.5Athan}. 
\end{proof}

\section{Hibi rings and flag varieties} \label{flags}
In this section, first we recall the definition and main properties of
Hibi rings associated with distributive lattices. Then we describe a
result of Reiner and Welker \cite{RW} that implies the existence of
Gorenstein initial complexes for ideals defining Gorenstein Hibi
rings. Finally we will illustrate Sagbi deformations of the coordinate
rings of flag varieties to certain Hibi rings.  For general facts on
Sagbi bases and Sagbi deformations we refer the reader to \cite{BC,
CHV} and \cite[Chapter 11]{Stu1}.

\subsection{Hibi rings} 

Let $(P,\leq)$ be a finite poset. If there is no danger of confusion
we will denote the poset only by the underlying set $P$.  A (possibly
empty) subset $I$ of $P$ is an order ideal if $x\leq y\in I$ implies
$x\in I$. The set $J(P)$ of the order ideals of $P$ is a poset under
set inclusion. This poset is a distributive lattice where join and
meet operations correspond to taking unions and intersections.  The
celebrated Birkhoff's theorem \cite{Birk} asserts that any finite
distributive lattice $L$ is lattice-isomorphic to $J(P)$ for some
poset $P$. Indeed one can take $P$ to be the set of {\em
join-irreducible} elements of $L$ with the poset structure induced by
$L$. An element $x\in L$ is called join-irreducible if it is not the
minimum of $L$ and cannot be written as $y\vee z$ for $z,y<x$. More
precisely, $J(P)\simeq L$ as lattices under the map sending any order
ideal $I=\{x_1,\dots,x_k\}$ of $P$ to $x_1\vee x_2\vee \dots \vee x_k$
where, by convention, the image of $\emptyset$ is the $\hat{0}$ of
$L$.

For any distributive lattice $L$ let $R_L$ be the polynomial ring over
the field $K$ whose variables are the elements of $L$.  For each pair
of incomparable elements $x,y\in L$ one defines the Hibi relation
$xy-(x\wedge y)(x\vee y)$.  The Hibi ideal $I_L$ is the ideal of $R_L$
generated by all the Hibi relations and the Hibi ring of $L$ is the
$K$-algebra defined by $I_L$:
$$I_L=( xy-(x\wedge y)(x\vee y) : x,y \mbox{ incomparable in } L)
    \quad \mbox{and} \quad H(L)=R_L/I_L.$$ Hibi proved in \cite{H}
    that $H(L)$ is a normal Cohen-Macaulay domain and is a homogeneous
    algebra with straightening law (ASL).  The main point of Hibi's
    proof is to describe $H(L)$ as a toric ring.  Let $Q$ be the {\em
    order polytope} of $P$, i.e., the convex hull of $\{ \chi_I : I
    \in J(P)\}$ where $\chi_I$ is the $0/1$ characteristic vector of
    $I$: $\chi_I(p) = 1$ if $p\in I$ and $\chi_I(p)=0$ otherwise.
    Birkhoff's theorem induces a $K$-algebra isomorphism
    $K[M(Q)]\simeq H(L)$. Here $K[M(Q)]$ is the monoid algebra
    associated to $Q$ as in Section \ref{intro}.  An important
    consequence of the ASL property is that $\ini_\succ(I_L)=(xy : x
    \mbox{ and } y \mbox{ incomparable})$ with respect to any reverse
    lexicographic order where $x \prec y$ whenever $x <_L y$. Hibi
    proved also that $H(L)$ is Gorenstein if and only if $P$ is graded
    i.e. all the maximal chains of $P$ have the same cardinality.

\begin{theorem}[Hibi] \label{hibithm} If $L$ is a distributive lattice then 
\begin{itemize}
\item[(1)]  the Hibi ring  $H(L)$ is a toric, normal, Cohen-Macaulay ASL,
\item[(2)]  the ideal $I_L$ has a quadratic squarefree initial ideal
  whose associated simplicial complex is the chain  complex of $L$,  and
\item[(3)] $H(L)$ is Gorenstein if and only if the poset of
  join-irreducible elements in $L$ is graded.
\end{itemize}
\end{theorem} 

Then Theorem~\ref{BRcor7} implies the following result.
  
\begin{theorem} \label{ReWe}
Let $L$ be a distributive lattice and assume that $H(L)$ is
Gorenstein. Then $I_L$ has a squarefree initial ideal which is
Gorenstein and whose associated simplicial complex is a cone over a
simplicial polytope.
\end{theorem}

A proof of Theorem \ref{ReWe} is given by Reiner and Welker in their
2002 preprint \cite{RW} .  We give a few details of their approach.
Let $P$ be the graded poset such that $L=J(P)$.  We may assume that
$P$ is a poset on $[n]$ and has rank $r$.  A chain of order ideals
$I_1 \subset I_2 \subset \cdots \subset I_t$ is called {\em
equatorial} if $f := \chi_{I_1} + \cdots + \chi_{I_t}$ has the
property that $\textup{min}_{p \in P} f(p) = 0$ and for every $j \in
[2,r]$, there exists a covering relation $p_{j-1} < p_j$ with
$p_{j-1}$ of rank $j-1$ and $p_j$ of rank $j$ such that $f(p_{j-1}) =
f(p_j)$. On the other hand, $I_1 \subset I_2 \subset \cdots \subset
I_t$ is {\em rank-constant} if $f$ is constant along ranks of $P$,
i.e., $f(p)=f(q)$ whenever $p$ and $q$ are elements of the same rank
in $P$.

\begin{definition} \cite[Definition 3.7]{RW} 
The equatorial complex $\Delta_{eq}(P)$ is the subcomplex of the order
complex $\Delta(J(P))$ whose faces are indexed by the equatorial
chains of non-empty order ideals in $P$.
\end{definition}

Theorem~3.6 in \cite{RW} proves that the collection of simplices $\{
\conv(\chi_I \,:\, I \in \mathcal R \cup \mathcal E) \}$ where
$\mathcal R$ (respectively $\mathcal E$) is a chain of non-empty rank
constant (equatorial) order ideals in $P$, gives a unimodular
triangulation of the order polytope $\mathcal O(P)$ called the
equatorial triangulation of $\mathcal O(P)$. Let $\mathcal O_{eq}(P)$
be the quotient polytope $\mathcal O(P) / V$ where $V$ is the linear
subspace spanned by the characteristic vectors of the rank-constant
ideals in $P$. This quotient polytope can be identified with the
orthogonal projection of $\mathcal O(P)$ onto $V^\perp$.  Reiner and
Welker show that the equatorial complex $\Delta_{eq}(P)$ can be
realized as the boundary complex of the simplicial polytope $Q$ that
is obtained by a reverse lexicographic triangulation of $\mathcal
O_{eq}(P)$ where the vertices of $\mathcal O_{eq}(P)$ corresponding to
order ideals $I$ with smaller cardinality come first.  This means the
following: if we order the vertices of $\mathcal O(P)$ reverse
lexicographically where those corresponding to the rank-constant order
ideals come first and then the rest is ordered according to the
cardinality of the order ideals, the unimodular triangulation we
obtain is the simplicial join of the simplex given by the
rank-constant ideals and $\Delta_{eq}(P)$. Now the initial ideal
$\ini_{eq}(I_L)$ of $I_L$ with respect to the above reverse
lexicographic term order is squarefree.  Moreover, the core of the
corresponding initial complex is $\Delta_{eq}(P)$ which is the
boundary complex of a simplical polytope.

Furthermore, Reiner and Welker show that $\ini_{eq}(I_L)$ is quadratic
if the width (i.e. the largest size of an antichain) of $P$ is at most
$2$ and need not be so if the width is larger than $2$.  In
particular, they give a positive answer to Question \ref{fque} for
Hibi ideals associated to graded posets $P$ of width at most $2$.
 
\subsection{  Flag varieties and their deformation to Hibi rings}

Let $V$ be a vector space of dimension $n$ over an algebraically
closed field $K$.  The Grassmann variety $G(m,n)$ is the set of
$m$-dimensional subspaces of $V$. It is a projective variety embedded
in the projective space $\PP^{N-1}$ where $N = {n \choose m}$ via the
Pl\"ucker map.  The coordinate ring $\Grass(m,n)$ of $G(m,n)$ in this
embedding is the $K$-subalgebra of $K[x_{ij} : 1\leq i\leq m, 1\leq
j\leq n]$ generated by the $m$-minors of the $m\times n$ matrix
$X=(x_{ij})$.  The algebra $\Grass(m,n)$ has a toric deformation to
the Hibi ring associated to the poset of maximal minors (see
\cite[Chapter 11]{Stu1}).  More generally, a similar statement holds
also for flag varieties. As we explain below, these toric deformations
are simple consequences of the straightening law for generic minors.
We first recall the definition of flag varieties and their
multi-homogeneous coordinate rings and then describe the toric
deformation in the language of Sagbi bases.
 
Consider a sequence $1\leq m_1<m_2<\dots<m_k<n$ and set
$M=\{m_1,\dots,m_k\}$. Define $F(M,n)=\{ V_1\subset V_2\subset
\dots\subset V_k \subset V : V_i \mbox{ a vector space of dimension }
m_i\}$.  Let $X=(x_{ij})$ be a $m_k\times n$ matrix of variables. For
$p\leq m_k$ and $a_1<\dots<a_p\leq n$ we denote by $[a_1,\dots,a_p]$
the $p$-minor of $X$ with row indices $1,2,\dots,p$ and column indices
$a_1,\dots,a_p$. If we set $L(M,n)=\{ [a_1,\dots,a_p] :
a_1<\dots<a_p\leq n,\,\, p\in M\}$, the multi-homogeneous coordinate
ring $\Flag(M,n)$ of $F(M,n)$ is
$$\Flag(M,n)=K[ [a_1,\dots,a_p] : [a_1,\dots,a_p] \in L(M,n)].$$

With the partial order $[a_1,\dots,a_p]\leq [b_1,\dots,b_q]$ if $p\geq
q$ and $a_i\leq b_i$ for $i=1,\dots,q$, the set of minors $L(M,n)$ becomes a
distributive lattice.
The straightening law for generic minors (see \cite{DRS} or \cite{BV})
asserts that the polynomial ring $K[x_{ij}]$ has a $K$-basis whose
elements are products of minors of $X$ of various order. It implies
immediately that $\Flag(M,n)$ has a $K$-basis $B(M,n)$ whose elements
are the products $\delta_1\dots \delta_v$ with $\delta_i\in L(M,n)$
and $\delta_1\leq \delta_2\leq \dots \leq \delta_v$.

Now comes the crucial (and easy) observation: if $\succ$ is a diagonal
term order and $\delta \neq \gamma\in B(M,n)$ then
$\ini_\succ(\delta)\neq \ini_\succ(\gamma)$. Recall that the {\em
initial algebra} of a $K$-algebra $R$ with respect to a term order
$\succ$, denoted as $\ini_\succ(R)$, is the $K$-vector space generated
by $\{\ini_\succ(f) \,:\, f \in R\}$. A subset $L$ of $R$ is a {\em
Sagbi basis} of $R$ with respect to $\succ$ if $\ini_\succ(R) =
K[\ini_\succ(\alpha) \,:\, \alpha \in L]$. The following result is 
part of the folklore in this subject.

\begin{proposition}  With the notation introduced above we have: 
\begin{itemize}
\item[(1)] the elements of $L(M,n)$ form a Sagbi basis of
  $\Flag(M,n)$, that is, the initial algebra
  $\ini_\succ(\Flag(M,n))=K[\ini_\succ(f) : f\in L(M,n)]$,
\item[(2)] the elements   $\ini_\succ(g)$ with $g\in B(M,n)$ form a
  $K$-basis   of  $\ini_\succ(\Flag(M,n))$, and 
\item[(3)] $\ini_\succ(\Flag(M,n))$ is the Hibi ring of $L(M,n)$. 
\end{itemize}
\end{proposition} 

\begin{proof} Let $f \neq 0 \in \Flag(M,n)$. Then $f$ can be
  written as a linear combination of elements in $B(M,n)$. Since
  distinct elements of $B(M,n)$ have distinct initial terms,
  $\ini_\succ(f)$ is the initial term of some element of $B(M,n)$. So
  any monomial in the initial algebra $\ini_\succ(\Flag(M,n))$ is of
  the form $\ini_\succ(g)$ for a unique $g\in B(M,n)$. This proves (1)
  and (2).  To prove (3) note that for $\delta,\gamma\in L(M,n)$ one
  has $\ini_\succ(\delta)\ini_\succ(\gamma)=\ini_\succ(\delta\wedge
  \gamma)\ini_\succ(\delta \vee \gamma)$. This gives a surjective
  $K$-algebra homomorphism from $H(L(M,n))$ to
  $\ini_\succ(\Flag(M,n))$. It must be also an isomorphism because the
  two rings have the same Hilbert function.
\end{proof}

So as we have seen, $\Flag(M,n)$ gets deformed to the Hibi ring
$H(L(M,n))$. One knows that $\Flag(M,n)$ is Gorenstein (even
factorial), see \cite[Chap.9]{Ful}.  This implies that the Hibi ring
$H(L(M,n))$ must be Gorenstein as well since it is a Cohen-Macaulay
graded domain with the Hilbert function of a Gorenstein ring
\cite[Cor.4.4.6]{BH}. One can also argue the other way around: check
that the poset of join-irreducible elements of $L(M,n)$ is graded
(indeed, it is a distributive lattice) and deduce that the Hibi ring
$H(L(M,n))$ is Gorenstein. Then, by Sagbi deformation, one has that
$\Flag(M,n)$ is Gorenstein. Below we present some examples describing
the poset of join -irreducible elements for some specific values of
$M$ and $n$.
 
Now we associate indeterminates $t_\alpha$ with $\alpha\in L(M,n)$ and
we obtain a presentation of $\Flag(M,n)$ as a quotient of $K[t_\alpha
: \alpha\in L(M,n)]$ via the map sending $t_\alpha$ to $\alpha$. The
kernel of this map is the Pl\"ucker ideal $\Plu(M,n)$:
$$\Flag(M,n)=K[t_\alpha : \alpha\in L(M,n)]/\Plu(M,n).$$
The generators of $\Plu(M,n)$ are quadrics which can be described in
terms of multilinear algebra, see \cite[Chap.9]{Ful}. In their reduced form  (in the sense of ASL theory or Sagbi basis theory)  they are of the form
$$t_\alpha t_\beta-t_{\alpha \vee \beta}t_{\alpha\wedge
    \beta}+\dots \mbox{ other terms  }  \lambda t_\gamma t_\delta $$
where
\begin{itemize}
\item[(1)] the $p$-minor $\alpha$ and  $q$-minor $\beta$ 
are incomparable in $L(M,n)$, and
\item[(2)] $\lambda\in \ZZ$ and in each term $\lambda t_\gamma
  t_\delta$ with $\lambda\neq 0$, $\gamma$ is a $p$-minor and
  $\delta$ is a $q$-minor with $\delta<\alpha\wedge \beta$ and
  $\gamma>\alpha\vee \beta$ and $\rank \alpha+\rank \beta= \rank
  \gamma+\rank \delta$.
\end{itemize} 

\begin{theorem}  \label{gorflag}
The ideal of Pl\"ucker relations $\Plu(M,n)$ defining
  the flag variety $F(M,n)$ has an initial ideal which is squarefree 
  and Gorenstein. 
\end{theorem}

\begin{proof} By Sagbi theory, any initial ideal of the
  ideal defining the initial algebra $H(L(M,n))$ is also an initial
  ideal of the ideal defining $\Flag(M,n)$.  So it is enough to show
  that the toric ideal $I_{L(M,n)}$ has a Gorenstein squarefree
  initial ideal. But this follows from Theorem~\ref{ReWe}.
\end{proof}

\begin{example} Consider the Grassmannian $\Grass(m,n) = \Flag(M,n)$
  with $M=\{m\}$. The join-irreducible elements of $L(M,n)$ are:
  $$\delta(a,b)=[1, 2,\dots,a-1,  a+b,a+1+b,\dots,m+b]$$ with
  $a=1,\dots,m$ and $b=1,\dots,n-m$. Note that $\delta(a,b)\leq
  \delta(c,d)$ if and only if $c\leq a$ and $b\leq d$ Hence the poset
  of join-irreducible elements $P$ of $L(M,n)$ is a $m\times (n-m)$
  grid, i.e.  the cartesian product of $[m]$ and $[n-m]$.  Therefore
  the width of $P$ is $\min(m,n-m)$.  It follows that Question
  \ref{fque} for $ \Grass(m,n)$ has a positive answer if
  $\min(m,n-m)\leq 2$. This is essentially the case $m=2$. 
  We analyze the case $m=2$ and $n=5$ in more detail. In this
  case $P = \{p<q\} \times \{1,2,3\}$ and  including the 
  empty order ideal there are ten order ideals of $P$ which
  we list using their maximal elements: 
  $$ \emptyset, \{p1\}, \{p2\}, \{p3\}, \{q1\}, \{q1,p2\}, \{q1,p3\}, \{q2\}, 
 \{q2,p3\}, \{q3\}.$$  
  We label these order ideals by $[12], [13], \ldots, [45]$ respectively.
  This ordering is consistent with the description of
  the join-irreducible elements $[13], [14], [15], [23], [34], [45]$.
  The order polytope is a six-dimensional polytope that is 
  the convex hull of the columns  of the matrix
$$\left[ \begin{array}{cccccccccc}
0 & 1 & 1 & 1 & 1 & 1 & 1 & 1 & 1 & 1 \\
0 & 0 & 1 & 1 & 0 & 1 & 1 & 1 & 1 & 1 \\
0 & 0 & 0 & 1 & 0 & 0 & 1 & 0 & 1 & 1 \\
0 & 0 & 0 & 0 & 1 & 1 & 1 & 1 & 1 & 1 \\
0 & 0 & 0 & 0 & 0 & 0 & 0 & 1 & 1 & 1 \\
0 & 0 & 0 & 0 & 0 & 0 & 0 & 0 & 0 & 1 
\end{array} \right], $$
where the columns correspond to the order ideals in the above order
and the rows correspond to $p1, p2, p3, q1, q2, q3$.  The Hibi
ideal $I_L$ is generated by five relations:
$$[14][23]-[13][24], [15][23]-[13][25], [15][24]-[14][25], [15][34]-[14][35], 
[25][34]-[24][35].$$ 
According to Theorem \ref{hibithm} the first terms of these
binomials are the minimal generators of the classical initial
ideal. The chain $[12] \subset [13] \subset [24] \subset [35] \subset [45]$
is the maximal chain of rank constant order ideals. And according
to Reiner-Welker construction we can take the following 
reverse lexicographic term order:
$$[12] \prec [13] \prec [24] \prec [35] \prec [45] \prec [14] \prec [23] \prec [15] \prec [25] \prec [34].$$
Then $\ini_\succ(I_L) = \langle [14][23], [14][25], [15][23], [15][34], [25][34]\rangle$, and the cone points of the corresponding simplicial
complex are precisely those in the maximal chain of rank constant
order ideals. Moreover, the core of this complex is the boundary 
complex of a pentagon whose
vertices are (in cyclic order) $[14], [34], [23], [25], [15]$.
%
\end{example}

\begin{example} Consider $\Flag(\{1,3,4,6\},7)$. The poset of 
join-irreducible elements of $L(M,n)$ is: 

 $$\begin{array}{cccccc} & & & & & 8\\ 234567 & 345678 & 4567 & 5678
   & 678 & 7\\ 134567 & 145678 & 1567 & 1678 & 178 & 1\\ 124567 &
   125678 & 1267 & 1278 & 128\\ 123567 & 123678 & 1237 & 1238 & 123
   \\ 123467 & 123478 & 1234 \\ 123457 & 123458
\end{array}
$$ 
Here $128$  stands for $[1,2,8]$ and so on. 
\end{example}

\section{Pfaffians} \label{jean-marie-pfaff}

In this section we let $X$ be an $n \times n$ skew symmetric matrix of
indeterminates: the diagonal entries of $X$ are zero, and $x_{ji} =
-x_{ij}$ for $i<j$. For $t = 2r \leq n$ and $J = \{1 \leq j_1 < \cdots
< j_t \leq n\}$ the $t$-minor of $X$ obtained from the columns and
rows indexed by $J$ is of the form $p_J(x)^2$ where $p_J(x)$ is a
polynomial of degree $r$. The polynomials $p_J(x)$ are called the
Pfaffians of order $t=2r$, and we let $\Pfaff(t,n)$ be the ideal
generated by all Pfaffians of order $t$. Squarefree initial ideals of
$\Pfaff(t,n)$ have been constructed \cite{HT}. However, even though
$\Pfaff(t,n)$ is always Gorenstein~\cite{Av} these initial ideals are
not. Here we give a sketch of the construction of Gorenstein initial
ideals from \cite{JW}. We thank Jakob Jonsson and Volkmar Welker for
generously sharing their manuscript in progress. We refer the reader
to this manuscript \cite{JW} for all the details.

The dimension and the degree of any initial ideal of
$\Pfaff(2(r+1),n)$ is equal to that of $\Pfaff(2(r+1),n)$: $\dim
K[X]/\Pfaff(2(r+1),n) =r(2n-2r-1)$ and
$$\deg K[X]/\Pfaff(2(r+1),n)= 
\prod_{1 \leq i \leq j \leq n-2r-1} \frac{2r+i+j}{i+j}.$$
The determinantal formula for the Hilbert series of $K[X]/\Pfaff(2(r+1),n)$
(see \cite{GK}) implies that the $a$-invariant is $-rn$. 

Suppose $p_J(x)$ is a Pfaffian of order $2r$ associated to the row and
column indices $J=\{1 \leq j_1 < \cdots < j_{2r} \leq n\}$.  Then the
terms of $p_J(x)$ are in bijection with the {\em perfect matchings} of
the complete graph on $2r$ vertices labeled by the elements of
$J$. For instance, if we take $n=6$, $r=2$, and $J= \{1, 2, 3, 4\}$,
the corresponding Pfaffian is $x_{14}x_{23} - x_{13}x_{24} +
x_{12}x_{34}$.  The terms of this Pfaffian correspond to the matchings
$\{1-4, \, 2-3\}$, $\{1-3, \, 2-4\}$, and $\{1-2, \, 3-4\}$,
respectively.  We will introduce a term order used in \cite{JW} that
picks as initial term that term of $p_J(x)$ corresponding to the
matching $\{j_1-j_{r+1}, \, j_2-j_{r+2}, \ldots, j_r-j_{2r}\}$.  For
this we let $d_{ij} = \min(j-i, n+i-j)$ for $i<j$. Now we totally
order the indeterminates so that $x_{ij} \prec x_{kl}$ whenever
$d_{ij} < d_{kl}$, and then we use a reverse lexicographic term order
induced by this ordering
\begin{example} \label{example7.1}
Let $n=6$ and $r=2$. We can use the following reverse 
lexicographic order:
$$\begin{array}{c}
x_{12} \prec x_{23} \prec x_{34} \prec x_{45} \prec x_{56} \prec x_{16} \prec\\
\ x_{13} \prec x_{24} \prec x_{35} \prec x_{46} \prec x_{15} \prec x_{26} \prec \\
x_{14} \prec x_{25} \prec x_{36}
\end{array} $$
The set of all Pfaffians is a Gr\"obner basis of $\Pfaff(4,6)$
where the underlined terms are the initial terms:
$$\begin{array}{c}
\underline{x_{36}x_{25}}-x_{26}x_{35}-x_{56}x_{23}, \, 
\underline{x_{36}x_{14}}-x_{46}x_{13}-x_{16}x_{34}, \, 
\underline{x_{25}x_{14}}-x_{15}x_{24}-x_{45}x_{12}, \\ 
\underline{x_{14}x_{26}}-x_{24}x_{16}-x_{46}x_{12}, \, 
\underline{x_{36}x_{15}}-x_{35}x_{16}-x_{13}x_{56}, \, 
\underline{x_{26}x_{15}}-x_{25}x_{16}-x_{56}x_{12}, \\  
\underline{x_{25}x_{46}}-x_{24}x_{56}-x_{26}x_{45}, \, 
\underline{x_{15}x_{46}}-x_{14}x_{56}-x_{16}x_{45}, \, 
\underline{x_{14}x_{35}}-x_{13}x_{45}-x_{15}x_{34}, \\ 
\underline{x_{46}x_{35}}-x_{36}x_{45}-x_{56}x_{34}, \, 
\underline{x_{36}x_{24}}-x_{26}x_{34}-x_{46}x_{23}, \, 
\underline{x_{35}x_{24}}-x_{25}x_{34}-x_{45}x_{23}, \\ 
\underline{x_{25}x_{13}}-x_{15}x_{23}-x_{35}x_{12}, \, 
\underline{x_{26}x_{13}}-x_{16}x_{23}-x_{36}x_{12}, \, 
\underline{x_{24}x_{13}}-x_{14}x_{23}-x_{34}x_{12} 
\end{array}
$$
\end{example}
\begin{proposition}[cf.  \cite{JW}] The initial term
of the Pfaffian $p_J(x)$ where $J = \{ j_1 < \cdots < j_{2r} \}$
is $x_{j_1j_{r+1}}x_{j_2j_{r+2}} \cdots x_{j_rj_{2r}}$. 
\end{proposition}
One consequence of this proposition is the following. Let $I(r,n)$ be
the squarefree ideal generated by the initial terms of the order $2r$
Pfaffians in $\Pfaff(2r,n)$.  This corresponds to a simplicial complex
$\Delta_{n, r-1}$ and the cone points of $\Delta_{n,r-1}$ correspond
to the variables which do not appear in the generators of $I(r,n)$.
These variables are of the form $x_{ij}$ where either $j-i \leq r-1$
or $n+i-j \leq r-1$. An easy counting argument shows that there are
$(r-1)n$ such cone points.  This is precisely
$-a(\Pfaff(2r,n))$. Therefore it is natural to ask whether $I(r,n)$ is
equal to $\ini_\succ(\Pfaff(2r,n))$.  The positive answer is the
content of the next result.

\begin{proposition} \cite[Theorem 2.1]{JW} The ideal $I(r,n)$ generated by
the initial terms of the Pfaffians in the ideal $\Pfaff(2r,n)$ is
equal to $\ini_\succ(\Pfaff(2r,n))$.
\end{proposition}

\begin{proof} Clearly $I(r,n) \subseteq \ini_\succ(\Pfaff(2r,n))$.
For the other inclusion we use Lemma \ref{kntrick}.  The simplicial
complex $\Delta_{n,r-1}$ corresponding to $I(r,n)$ can be described as
follows: Let $\Omega_{n} = \{(i,j)\, : \, 1 \leq i < j \leq n\}$ which
we will think of as the edges and diagonals of a convex $n$-gon. For
$j \geq 1$, a $j$-crossing is a subset of $j$ elements of $\Omega_n$
which mutually intersect and where all $2j$ endpoints are
distinct. Then $\Delta_{n,r-1}$ is the simplicial complex of all
subsets of $\Omega_n$ which do not contain an $r$-crossing. Observe
that the minimal nonfaces of $\Delta_{n,r-1}$ are precisely
$r$-crossings, and they correspond to the minimal generators of
$I(r,n)$.  By the results in \cite{DKM} and \cite{Jo} the simplicial
complex $\Delta_{n,r}$ is a pure complex of dimension $r(2n-2r-1)-1$
and has $\prod_{1 \leq i \leq j \leq n-2r-1} \frac{2r+i+j}{i+j}$
facets.  Now Lemma \ref{kntrick} implies that $I(r+1, n) =
\ini_\succ(\Pfaff(2(r+1),n))$.
\end{proof}
Further results in \cite{DGJM} show that the core $\Delta_{n,r}'$ of
$\Delta_{n,r}$ is a simplicial sphere. With this result 
we get the main theorem of this section.
\begin{theorem} \cite[Theorem 2.1]{JW} The ideal $\ini_\succ(\Pfaff(2r,n))$ is 
a squarefree Gorenstein  initial ideal. 
\end{theorem}
We finish this section by pointing out that $\Delta_{n,1}'$ is 
the boundary complex of the $n$-associahedron, and hence it
is a polytopal sphere \cite{DGJM}. It remains open whether
$\Delta_{n,r}'$ is a polytopal sphere in general. 

\section{Minors} \label{minors}
In this section, we return to Question~\ref{q1} and illustrate a
family of determinantal ideals such that for each $I$ in the family
there is an initial ideal $\ini_\succ(I)$ with the same Betti numbers
as $I$.

\begin{theorem}\label{n-1 in nxn} Let $I$ be the ideal of $(n-1)$-minors 
  of the generic $n\times n$ matrix $X=(x_{ij})$ with $n>2$. Set $V=\{
  x_{ij} : 0\leq j-i\leq 1\}\cup \{x_{n1}\}$ and $W = \{x_{ij} :
  x_{ij}\not\in V\}$. Let $Y$ be the matrix obtained from $X$ by
  replacing $x_{ij}$ with $0$ if $x_{ij} \in W$. Let $\succ$ be any
  reverse lexicographic order on the $x_{ij}$ 
  such that $x_{ij}>x_{hk}$ if $x_{ij}\in
  V$ and $x_{hk}\in W$. Then $\ini_\succ(I)$ is a square-free monomial
  ideal with Betti numbers equal to those of $I$ and the core of the
  associated initial complex is the cyclic polytope with $2n$ vertices
  in $\RR^{2n-4}$. More precisely, $\ini_\succ(I)$ is the
  specialization of $I$ by the regular sequence $W$ or in other words,
  $\ini_\tau(I)$ is the ideal of $(n-1)$-minors of $Y$.
\end{theorem} 

We note that part of Daniel Soll's thesis \cite{Soll} has results 
about the initial complexes of determinantal ideals, and a result
similar to Theorem \ref{n-1 in nxn} appears there as well. 

Given a matrix $Z$ we define a graph $G(Z)$ as follows. The vertices of
$G(Z)$ are the elements $z_{ij}$ such that $z_{ij}\neq 0$ and the
edges are the pairs $\{z_{ij},z_{hk}\}$ such that $i=h$ or $j=k$. Note
that $G(Y)$ is a cycle of length $2n$. The statement of
Theorem~\ref{n-1 in nxn} remains true whenever $V$ is a subset of the
$x_{ij}$'s such that the corresponding matrix $Y$ has the property
that the graph $G(Y)$ is a cycle of length $2n$. The main ingredient
needed in the proof of Theorem \ref{n-1 in nxn} is the following
lemma.

\begin{lemma}\label{main-1} Suppose $Y$ is a $n\times n$ matrix such
  that $G(Y)$ is a cycle of length $2n$.  Let $J_k$ be the ideal of
  $k$-minors of $Y$. Then for all $k<n$,
  $$J_k=( \prod_{v\in A} v : A \mbox{ is an independent set of } G(Y) \mbox{
    with } |A|=k).$$
\end{lemma} 

\begin{proof} [Proof of Theorem~\ref{n-1 in nxn}] With the notation of
Theorem~\ref{n-1 in nxn}, let $J$ be the ideal generated by the
$(n-1)$-minors of $Y$.  From Lemma~\ref{main-1}, $J$ is a square-free
monomial ideal. Let $\Delta$ be the initial complex of $J$.
Lemma~\ref{main-1} implies that the facets of $\core(\Delta)$ are the
sets obtained as unions of $n-2$ disjoint edges of $G$.  Since the
size of any such facet is $2n-4$ and $Y$ has $2n$ nonzero entries,
the codimension of $J$ is $4$. Now using the facts
that $I$ defines a Cohen-Macaulay ring and that $I+(W)=J+(W)$, we 
conclude that $W$ is
a regular sequence modulo $I$.  This in turn shows that $J$ is a
specialization of $I$ by the regular sequence $W$ and hence $J$ and
$I$ have the same Betti numbers and same Hilbert function. Since $J
\subset \ini_\succ(I)$ holds by construction it must then be that $J =
\ini_\succ(I)$. That $\core(\Delta)$ is the cyclic polytope with
$2n$ vertices in $\RR^{2n-4}$ follows from the facet description given
above and Gale's evenness characterization of the facets of the cyclic
polytope; see \cite[Chapter 0]{Z}. 
\end{proof}

It remains to prove Lemma~\ref{main-1}. To this end, let us introduce
some notation.  Let $Z$ be a matrix such that each row and column of
$Z$ contains at most two non-zero entries.  Then each vertex of $G(Z)$
is contained in at most two edges. Therefore the connected components
of $G(Z)$ are either paths or cycles. The decomposition of $G(Z)$ into
connected components correspond to a block decomposition of $Z$ as
follows: if $G(Z)$ has connected components $G_1,\dots,G_r$ then, up
to row and column permutations and after eliminating zero rows and
columns from $Z$, $Z$ has a block decomposition of the form:

$$ 
\begin{array}{cccccc}
Z_1   & 0 &     \dots  & 0  \\
0     & Z_2 & 0 &  \dots  \\ 
  &  \ddots    & \ddots     &  \ddots \\
0  &\dots       &   0     &  Z_r \\
\end{array}
$$  

\begin{itemize}
\item[(1)] If $G_i$ is a cycle, then it is a cycle of even length with
  vertices $y_1,\dots,y_{2k}$, and $Z_i$ is the $k\times k$
  matrix
$$
\begin{array}{cccccc}
y_1   & y_2 &  0   & \dots & 0    \\
0     & y_3 & y_4 & 0& \dots  \\ 
  &      & \ddots     &  \ddots \\
0      & \dots    &   0  &  y_{2k-3} & y_{2k-2}   \\
y_{2k} &  0   & \dots    &    0       & y_{2k-1}  
\end{array}
$$

\item[(2)] If $G_i$ is a path of odd length with vertices
  $y_1,\dots,y_{2k-1}$, then $Z_i$ is the $k\times k$ matrix
$$
\begin{array}{cccccc}
y_1   & y_2 &  0   & \dots    & 0  \\
0     & y_3 & y_4 & 0& \dots  \\ 
  &      & \ddots     &  \ddots \\
0      & \dots    &   0  &  y_{2k-3} & y_{2k-2}   \\
0      &  0   & \dots    &    0       & y_{2k-1}  
\end{array}
$$
or its transpose. 

\item[(3)] If $G_i$ is a path of even length with vertices
  $y_1,\dots,y_{2k}$ then $Z_i$ is the $k\times (k+1)$ matrix
$$
\begin{array}{cccccc}
y_1   & y_2 &  0   & \dots & & 0  \\
0     & y_3 & y_4 & 0& \dots  \\ 
  &      & \ddots     &  \ddots \\
0      & \dots    &   0  &  y_{2k-3} & y_{2k-2}   &0  \\
0      &  0   & \dots    &    0       & y_{2k-1}   & y_{2k}
\end{array}
$$
or its transpose. 
\end{itemize} 

It follows that if $Z$ is a square matrix containing no zero rows or
columns, then $\det Z=0$ if one of the $G_i$ is a path of even length.
Otherwise, $\det Z$ is (up to sign) the product of the determinants
$\det Z_i$ associated to the blocks. Furthermore, $\det Z_i=y_1y_3\dots
y_{2k-1}-y_2y_4\dots y_{2k}$ if $G_i$ is a cycle (case (1) above) and
$\det Z_i=y_1y_3\dots y_{2k-1}$ if $G_i$ is a path of odd length (case
(2) above).

\begin{proof}[Proof of \ref{main-1}]  Set $G=G(Y)$. Denote by $V$ the
  set of vertices of $G$. For simplicity, we identify square-free
  monomials in the variables in $V$ with subsets of $V$. Denote by
  $U_k$ the ideal generated by the independent subsets of $G$ of
  cardinality $k$. We have to show that $J_k=U_k$ for all $k<n$.  
  
  The inclusion $J_k\subseteq U_k$ follows from the very definition of
  determinant.  For the other inclusion note that if $Z$ is the
  $k\times k$ sub-matrix of $Y$ with row indices $R =
  \{r_1,\dots,r_k\}$ and column indices $C=\{c_1,\dots,c_k\}$ then
  $G(Z)$ is the subgraph of $G$ whose vertices $y_{ij}$ satisfy $i\in
  R$ and $j\in C$.  In particular, if $k<n$ then $G(Z)$ is not a cycle
  and so its connected components are lines. It follow that if $k<n$
  then $\det Z$ is either $0$ or a monomial in the variables of $V$.
  Consider now an independent set of cardinality $k<n$ of $G$, say
  $y_{i_1j_1},\dots, y_{i_kj_k}$.  By construction, $i_a\neq i_b$ and
  $j_a\neq j_b$ if $a\neq b$.  Consider the sub-matrix $Z$ of $Y$ with
  row indices $i_1,\dots,i_k$ and column indices $j_1,\dots,j_k$.  By
  construction $y_{i_1j_1}\cdots y_{i_kj_k}$ appears in $\det Z$ and,
  since we know that $\det Z$ is either $0$ or a monomial, we may
  conclude that $\det Z$ is $y_{i_1j_1}\cdots y_{i_kj_k}$ up to sign.
  This implies that $U_k\subseteq J_k$ and concludes the proof.
\end{proof}


\begin{thebibliography}{99}

\bibitem{A} C.Athanasiadis. Ehrhart polynomials, simplicial polytopes,
magic squares and a conjecture of Stanley, {\em J. Reine
Angew. Math.} {\bf 583} (2005), 163-174.

\bibitem{Av} L.Avramov. A class of factorial domains, {\em Serdica}
{\bf 5} (1979), no. 4, 378--379.

\bibitem{Birk} G.Birkhoff.  {\em Lattice theory}, American Mathematical
  Society Colloquium Publications, Vol. XXV, American Mathematical
  Society, Providence, R.I. 1967. 

\bibitem{BLSWZ} A.Bj\"orner, M.Las Vergnas, B.Sturmfels, N.White and 
G.Ziegler. {\em Oriented matroids}. Second edition. Encyclopedia
of Mathematics and its Applications, 46. Cambridge University Press, 
Cambridge, 1999.
  
\bibitem{BC} W.Bruns and A.Conca. Gr\"obner bases and determinantal
    ideals, {\em Commutative algebra, singularities and computer
    algebra (Sinaia, 2002)}, 9--66, NATO Sci. Ser. II Math.  Phys.
    Chem., 115, Kluwer Acad. Publ., Dordrecht, 2003.

\bibitem{BR} W.Bruns and T.R\"omer. $h$-vectors of Gorenstein
    polytopes, math.AC/0508392.

\bibitem{BH} W. Bruns and J. Herzog. {\em Cohen-Macaulay Rings},
Cambridge University Press, Cambridge, 1996.

\bibitem{BV} W.Bruns and U.Vetter. { \em Determinantal Rings}, Lecture
Notes in Mathematics  {\bf 1327}, Springer-Verlag, Berlin, 1988.

\bibitem{C} A.Conca. Gr\"obner bases of ideals of minors of a
symmetric matrix,  {\em J. Algebra} {\bf 166} (1994), 406-421.

\bibitem{CHV} A.Conca, J.Herzog and G.Valla. Sagbi bases with
    applications to blow-up algebras, {\em J. Reine Angew. Math.}
  {\bf 474} (1996), 113--138.

\bibitem{Cocoa} CoCoA Team. {\em CoCoA: a system for doing
     Computations in Commutative Algebra}, Available at
     http://cocoa.dima.unige.it.
 
\bibitem{DRS} P.Doubilet, G.C.Rota and J.Stein. {\em On the foundations of
    combinatorial theory IX},  Combinatorial methods in Invariant
    Theory. Studies in Appl. Math. 53 (1974), 185--216. 
    
\bibitem{DGJM} A.Dress, S.Gr\"unewald, J.Jonsson and V.Moulton. 
Paper in preparation.

\bibitem{DKM} A.Dress, J.Koolen and V.Moulton. $4n-10$, {\em
Annals of Combinatorics} {\bf 8} (2005) 463--471.

\bibitem{E} D.Eisenbud. {\em Commutative Algebra with a view
       towards Algebraic Geometry}, Graduate Texts in Mathematics
     {\bf 150}, Springer-Verlag, New York, 1995.

\bibitem{HT} J.Herzog and N.V.Trung. Gr\"obner bases and
      multiplicity of determinantal and Pfaffian ideals, {\em Adv.
      Math.} {\bf 96} (1992), no. 1, 1--37.

\bibitem{Ful} W.Fulton. {\em Young tableaux,  with applications to
    representation theory and geometry},  London Mathematical Society
    Student Texts {\bf 35}, Cambridge University Press, Cambridge, 1997.

\bibitem{GK} S.R.Ghorpade and C.Krattenthaler. The Hilbert
Series of Pfaffian rings, {\em Algebra, arithmetic and geometry with 
applications} (West Lafayette, IN, 2000), 337--356, Springer, Berlin, 
2004. 

\bibitem{Goto} S.Goto. On the Gorensteinness of determinantal loci,
{\em J. Math. Kyoto Univ.} {\bf 19} (1979), no. 2, 371--374.

\bibitem{GW} S.Goto and K.Watanabe. On graded rings. {\em
  I. J. Math. Soc. Japan} {\bf 30} (1978), no. 2, 179--213.

\bibitem{H} T.Hibi. Distributive lattices, affine semigroup rings and
algebras with straightening laws, {\em Commutative Algebra and
Combinatorics} (Kyoto, 1985), 93--109, Adv. Stud. Pure Math., {\bf
11}, North-Holland, Amsterdam, 1987.

\bibitem{HO2005} T. Hibi and H. Ohsugi. Special simplices 
and Gorenstein toric rings, math.AC/0503666.
  
\bibitem{J} A. Jensen. {\em Gfan, a software system for Gr\"obner fans},
available at {\tt http://home.imf.au.dk/ajensen/software/gfan/gfan.html}.

\bibitem{Jo} J.Jonsson. Generalized triangulations and diagonal-free 
subsets of stack polyominoes,  {\em J. Combin. Theory Ser. A} 
{\bf 112} (2005), 117--142.

\bibitem{JW} J.Jonsson and V.Welker. A spherical initial ideal
of Pfaffians, Manuscript (2005). 

\bibitem{KL} H.Kleppe and D.Laksov. The algebraic structure and
deformation of Pfaffian schemes, {\em J. Algebra} {\bf 64} (1980), 
no. 1, 167--189.

\bibitem{RW} V.Reiner and V.Welker. On the Charney-Davis and
Neggers-Stanley conjectures,  {\em J. Combin. Theory Ser. A}  {\bf 109}  
(2005),  247--280. 


\bibitem{Soll} D.Soll. Diploma Thesis, in preparation. 

\bibitem{Stan80} R. Stanley. Decompositions of rational convex polytopes,
{\em Annals of Discrete Math.} {\bf 6} (1980), 333-342.  

\bibitem{Sta} R.Stanley. {\em Combinatorics and Commutative Algebra}, 
Progress in Mathematics {\bf 41}, Birkh\"auser, Boston 1996.  

\bibitem{Stu} B.Sturmfels. Gr\"obner bases and Stanley
decompositions of determinantal rings, {\em Math.Z.} {\bf 205} (1990), 
137-144. 

\bibitem{Stu1} B.Sturmfels. {\em Gr\"obner Bases and Convex
Polytopes}, American Mathematical Society, Providence, RI, 1995.

\bibitem{Stu2} B.Sturmfels, {\em Algorithms in Invariant Theory},
Springer Verlag, Vienna, 1993.

\bibitem{Vil} R.Villarreal. {\em Monomial algebras}, Monographs and
Textbooks in Pure and Applied Mathematics, 238. Marcel Dekker Inc., 
New York, 2001. 

\bibitem{Z} G.Ziegler, {\em Lectures on Polytopes}, Graduate Texts in
Mathematics, Springer-Verlag, Berlin, 1995.
 
\end{thebibliography}
\end{document}